\documentclass[12pt]{article}
\usepackage{amsfonts,latexsym,color}
\usepackage{latexsym}
\usepackage{srcltx}
\usepackage{amsmath, amscd, amssymb}
\begin{document}
\baselineskip 16pt

\def \Aut {\mbox{\rm Aut\,}}
\def \Sym {\mbox{\rm Sym\,}}
\def \Stab {\mbox{\rm Stab\,}}
\def \Rot {\mbox{\rm Rot\,}}
\def \Fix {\mbox{\rm Fix\,}}
\def \Mon {\mbox{\rm Mon\,}}
\newfont{\sBbb}{msbm7 scaled\magstephalf}
\newcommand{\BC}{\mbox{\Bbb C}}
\newcommand{\BR}{\mbox{\Bbb R}}
\newcommand{\BZ}{\mbox{\Bbb Z}}
\newcommand{\sZ}{\mbox{\sBbb Z}}

\newcommand{\Iso}{\mbox{Iso}\,}
\newcommand{\Fi}{\mbox{Fix}\,}
\renewcommand{\labelenumi}{\theenumi}
\newcommand{\qed}{\mbox{\raisebox{0.7ex}{\fbox{}}}}
\newtheorem{theorem}{Theorem}[section]
\newtheorem{problem}{Problem}
\newtheorem{defin}{Definition}
\newtheorem{lemma}[theorem]{Lemma}
\newtheorem{conj}{Conjecture}
\newtheorem{prop}[theorem]{Proposition}
\newtheorem{op}{Open Problem}
\newtheorem{example}{Example}[section]
\newtheorem{note}{Note}
\newtheorem{remark}{Remark}
\newtheorem{corollary}[theorem]{Corollary}
\newenvironment{pf}{\medskip\noindent{Proof:}
  \hspace{-.5cm}      \enspace}{\hfill \qed \newline \smallskip}
\newenvironment{pflike}[1]{\medskip\noindent {#1}
   \enspace}{\medskip}
\newenvironment{defn}{\begin{defin} \em}{\end{defin}}
\newenvironment{nt}{\begin{note} \em}{\end{note}}
\newenvironment{rem}{\begin{remark} \em}{\end{remark}}
\newenvironment{examp}{\begin{example} \rm}{\end{example}}
\newcommand{\lra}{\longrightarrow}
\newcommand{\vect}[2]{\mbox{$({#1}_1,\ldots,{#1}_{#2})$}}
\newcommand{\comb}[2]{\mbox{$ \left( \begin{array}{c}
        {#1} \\ {#2} \end{array}\right)$}}
\newcommand{\tve}[1]{\mbox{$\mathbf{\tilde {#1}}$}}
\newcommand{\bve}[1]{\mbox{$\mathbf{{#1}}$}}
\setlength{\unitlength}{12pt}
\renewcommand{\labelenumi}{(\theenumi)}
\def\mod{\hbox{\rm mod }}

\title{
Classification of  nonorientable regular embeddings of complete
bipartite graphs }

\author{Jin Ho Kwak  \\
{\small Mathematics, Pohang University of Science and Technology,
Pohang,
790-784 Korea}\vspace{3mm} \\
Young Soo Kwon  \\
 {\small Mathematics, Yeungnam University, Kyeongsan, 712-749 Korea
} }
\date{}
\maketitle

\renewcommand{\thefootnote}{\empty}
\footnote{The first and the second authors are supported by the Korea Research
Foundation Grant funded Korean Government KRF-2007-313-C00011 and KRF-2008-331-C00049, respectively.}

\vspace{-.6cm}

\begin{abstract}
A  2-cell embedding of a graph $G$ into a closed (orientable or
nonorientable) surface is called regular if its automorphism group
acts regularly on the flags - mutually incident vertex-edge-face
triples. In this paper, we classify the regular embeddings of
complete bipartite graphs $K_{n,n}$ into nonorientable surfaces. Such
a regular embedding of $K_{n,n}$ exists only when $n =
2p_1^{a_1}p_2^{a_2}\cdots p_k^{a_k}$ (a prime decomposition of
$n$) and all $p_i \equiv
\pm 1 (\mod 8)$. In this case, the
number of those regular embeddings of $K_{n,n}$ up to isomorphism
is $2^k$.   \vskip10pt \noindent{\bf Keywords:} Graph,
surface, regular embedding, regular map 
\end{abstract}

\section{Introduction}

A 2-cell embedding of a connected graph into a  connected closed
surface is called a topological {\it map}. An automorphism of a
map is an automorphism of the underlying graph which can be
extended to a self-homeomorphism of the supporting surface in the
embedding. The automorphisms of a map $\cal M$ act semi-regularly
on its {\em flags} - mutually incident vertex-edge-face triples.
If the automorphism group $\Aut({\cal M})$ of
$\cal M$ acts regularly on the flags then the map $\cal M$ as well as
the corresponding embedding are also called {\it regular}. For a given
map $\cal M$ whose supporting surface is orientable, the set
$\Aut^{+}({\cal M})$ of orientation-preserving automorphisms of
the map $\cal M$ acts semi-regularly on its {\em arcs} - mutually
incident vertex-edge pairs. If  $|\Aut^{+}({\cal M})|$ acts regularly on its arcs then the map $\cal M$ is
called {\it orientably regular}. Therefore, if a supporting
surface is orientable, a regular map means an orientably regular
map having an orientation-reversing automorphism.

 There is a combinatorial
description of topological  maps. For a given finite set $F$ and
three fixed-point free involutary permutations   $\lambda, \rho$ and
$\tau$ on $F$, a quadruple ${\cal M} = (F;\,\lambda,\rho,\tau)$ is
called a \emph{combinatorial map} if they satisfy the following
properties:
\begin{enumerate}
  \item [{\rm (i)}] $\lambda\tau=\tau\lambda $ \quad and
  \vspace{-2mm}
  \item [{\rm (ii)}] the group $\langle\lambda,\rho,\tau\rangle$ generated by $\lambda, \rho$ and $\tau$ acts transitively
on $F$.
\end{enumerate}
A set $F$ is called the \emph{flag set}, and $\lambda, \rho$ and
$\tau$ are called the {\it longitudinal\/},  {\it rotary\/} and {\it
transversal involutions\/}, respectively.  The group
$\langle\lambda,\rho,\tau\rangle$ is called the {\it monodromy
group\/} of ${\cal M}$ and denoted by Mon$({\cal M})$.

 For a
given combinatorial map ${\cal M} = (F;\,\lambda,\rho,\tau)$, a
corresponding topological map is constructed as follows. The
orbits under the subgroups $\langle\rho,\tau\rangle$,
$\langle\lambda,\tau\rangle$ and $\langle\rho,\lambda\rangle$ are
vertices,  edges and face-boundaries of  $\cal M$, respectively.
The incidence is defined by a nontrivial set intersection. The
flag set $F$ in a combinatorial map corresponds to the set of
mutually incident vertex-edge-face triples in a topological map.
In fact, every topological map can be described by a combinatorial
map and vice versa. In this paper, when ${\cal M}=(F;\,\lambda,\rho,\tau)$ is a combinatorial description of an embedding ${\cal M}_1$ of $G$, we assume that  $F$ is the flags set of ${\cal M}_1$ and  the underlying graph of ${\cal M}$ is $G$.  For a detail, the reader is
referred to the paper \cite{N}.

A combinatorial map ${\cal M}=(F;\,\lambda,\rho,\tau)$ is said to be
{\it nonorientable} if
 the even-word subgroup
$\langle\rho\tau,\tau\lambda\rangle$ of Mon$({\cal M})$ acts
transitively on $F$; otherwise it is said to be
 {\it orientable\/}. In fact, a combinatorial map
 is nonorientable if and only if the supporting surface of the corresponding topological map
 is nonorientable.

The homomorphisms, isomorphisms and automorphisms of combinatorial
maps are described as follows. For given two combinatorial maps
${\cal M}=(F;\,\lambda,\rho,\tau)$ and ${\cal M}'=(F';\,\lambda
',\rho ',\tau ')$ , a {\it map homomorphism} $\psi: {\cal M}\to
{\cal M}'$  is a function $\psi:F\to F'$ satisfying
$$\psi\lambda =\lambda '\psi,\quad\psi\rho=\rho '\psi\quad\text{ and }
\quad\psi\tau=\tau '\psi.$$
Since a monodromy group is assumed to
act transitively on the flag set, the function $\psi$ is
surjective. If it is one-to-one then it is called a {\it map
isomorphism}. Furthermore, if ${\cal M} = {\cal M}'$ then an
isomorphism of the map is called a {\it map automorphism}. Let
$\Aut({\cal M})$ denote the set of  map automorphisms of ${\cal
M}$. Then, $\Aut({\cal M})$ is nothing but the centralizer of
$\Mon({\cal M})$ in $S_F$, the symmetric group on $F$. It means
that $\Aut({\cal M})$ is a group under composition.
 Since  the monodromy
group Mon$({\cal M})$ acts transitively on $F$,  $\Aut({\cal M})$
acts semi-regularly on  $F$. So, we have
  $|\Mon({\cal M})|\ge |F|  \geq |\Aut({\cal
M})|$ for every combinatorial map ${\cal M}$. If one of the
equalities $|\Mon({\cal M})| = |F|$ and $|F| = |\Aut({\cal M})|$
holds then so does the other. In this case,
 both Mon$({\cal M})$ and $\Aut({\cal
M})$ act regularly on $F$ and the  map ${\cal M}$ is  said to be
{\it regular\/}.
 As a well-known result in permutation group
theory (see Theorem 6.5 in \cite{Hu}), one can see that if $\cal{M}$
is regular then the  associated permutation groups
 $\Aut (\cal{M})$ and $\Mon (\cal{M})$ are isomorphic.
 The isomorphism can be understood with a fixed flag $\xi$
 as follows: Since $\Aut({\cal M})$ acts regularly on $F$, the
flag set $F$ can be identified with $\Aut({\cal M})$ as a set such
that $\xi$ is identified with the identity automorphism. Then, two
actions of $\Aut({\cal M})$ and Mon$({\cal M})$ on $F$ are
equivalent to the left regular and the right regular representations
of $\Aut({\cal M})$, respectively.

One of the standard problems in topological graph theory is a
classification of  orientably regular embeddings or regular
embeddings of a given class of graphs. In recent years, there has
been particular interest in the orientably regular embeddings and
regular embeddings of complete bipartite graphs $K_{n,n}$ by
several authors \cite{DJKNS1, DJKNS2,J, JNS1, JNS2, KK, KK1, NSZ}.
The regular (or reflexible regular) embeddings and self-Petrie
dual regular embeddings  of $K_{n,n}$ into orientable surfaces
were classified by the authors \cite{KK}. During preparing this
paper, G. Jones ~\cite{J} informed us that the classification of
orientably regular embeddings of $K_{n,n}$ has been completed. In
this paper, we classify the nonorientable regular embeddings of
$K_{n,n}$.

The nonorientable regular embeddings of complete graphs have been
classified by S. Wilson \cite{Wi}: There exists a nonorientable
regular embedding of $K_n$ exists only when $n=3,4$ or $6$. For
the $n$-dimensional cube $Q_n$, the classification of
nonorientable regular embeddings has been done
 by R. Nedela and the second
author \cite{KN} by showing a nonexistence of a nonorientable
regular embedding of $Q_n$ except  $n =2$. Contrary to all known
cases, there exists a nonorientable regular embedding of $K_{n,n}$
for infinitely many $n$.   The following theorem is the main
result in this paper.

\begin{theorem} \label{main-theorem}
For any integer $n \equiv 0, 1$ or $3 (\mod 4)$,  no nonorientable
regular embedding of $K_{n,n}$ exists. For $n =
2p_1^{a_1}p_2^{a_2}\cdots p_k^{a_k}$ (the prime decomposition of
$n$), the number of nonorientable regular embeddings of $K_{n,n}$
up to isomorphism is  $2^k$ if $p_i \equiv \pm 1 (\mod 8)$ for all
$i=1,2,\ldots,k$; $0$ otherwise.
\end{theorem}

This paper is organized as follows. In the next section, we
introduce a triple of three graph automorphisms of $G$, called an
{\em admissible triple} for $G$, corresponding to $\lambda, \rho,
\tau$ for a regular embedding of $G$. In Section $3$, we construct
some nonorientable regular embeddings of $K_{n,n}$ by forming
admissible triples for $K_{n,n}$. In the last two sections,
Theorem~\ref{main-theorem} is proved by showing that no other
nonorientable regular embedding of $K_{n,n}$ exists beyond those
constructed in Section 3 up to isomorphism.


\section{ An admissible triple of graph automorphisms}

 For a given
graph $G$, the automorphism group $\Aut(G)$ acts faithfully on
both the vertex set $V(G)$ and the arc set $D(G)$. Moreover, in an
embedding ${\cal M}$ of $G$ whose vertex valencies are greater than two,
$\Aut({\cal M}) \le \Aut(G)$ and $\Aut({\cal M})$ acts faithfully
on the flag set $F({\cal M})$. So, we consider a graph
automorphism as a permutation
 of $V(G)$, $D(G)$ or $F({\cal M})$ according to the context.

Let $G$ be a graph and let ${\cal M}=(F;\,\lambda,\rho,\tau)$ be a
combinatorial description of a regular embedding of $G$. For our
convenience, we fix an incident vertex-edge pair $(v,e)$ which is
extended to a flag $\xi=(v,e,f),$  called a root. Then, there
exist three involutory graph automorphisms $\ell, r$ and $t$ in
$\Aut(G)$ which correspond to the three involutions $\lambda,
\rho$ and $\tau$ in $\Mon({\cal M})$ with a root flag $\xi$. It
means that $\ell(\xi)=\lambda\cdot\xi,\ r(\xi)=\rho\cdot\xi \ \
\mbox{and} \ \ t(\xi)=\tau\cdot\xi$. Note that the three
involutory graph automorphisms $\ell, r$ and $t$ of $G$ satisfy
the following properties.
\begin{itemize}
\item [{\rm (i)}] $\Gamma = \langle \ell, r,t \rangle$ acts
transitively  on the arc set $D(G)$.\vspace{-.3cm}
  \item  [{\rm (ii)}] The stabilizer $\Gamma_v$ of the vertex $v$ is
$\Gamma_v=\langle r,t \rangle$ and it is isomorphic to the
dihedral group $D_n$,
 and its cyclic subgroup $\langle r t \rangle$ acts regularly on
 the
arcs emanating from $v$, where $n$ is the valency of
$G$.\vspace{-.3cm}
 \item  [{\rm (iii)}] The stabilizer $\Gamma_e$ of the edge $e$ is
$\Gamma_e=\langle \ell,t \rangle$ and it is isomorphic to
 the Klein four-group $\mathbb{Z}_2 \times
\mathbb{Z}_2$.
\end{itemize}
We call such a triple $( \ell, r,t)$ of involutory automorphisms
of $G$ an {\it admissible triple for  a regular embedding of $G$}
or simply an {\it admissible triple for $G$}. The vertex $v$ and
the edge $e$ are called
 a {\em root vertex} and a
{\em root edge}, respectively. Note that $\Aut({\cal M}) = \langle
\ell, r, t \rangle$.

Conversely, for a given admissible triple $(\ell, r, t)$ for $G$
with the root vertex $v$, one can construct a combinatorial map
${\cal M}_1=(F_1;\,\lambda_1,\rho_1,\tau_1)$ as follows: Let $F_1
= \langle \ell, r, t \rangle$ as a set and  for any $g \in F_1$,
let
$$ \lambda_1\cdot g =g\ell, \ \ \rho_1\cdot g =g r
\ \ \mbox{and} \ \ \tau_1\cdot g=g t,$$ namely, the right
translations by  $\ell, r, t$, respectively.  By conditions (i)
and (ii) above, one can show that $|F_1| = 4|E(G)|$ and
$\Mon({\cal M}_1) =\langle\lambda_1,\rho_1,\tau_1 \rangle$ acts
regularly on $F_1$. Consequently, the map ${\cal M}_1$ becomes a
regular map. We denote this regular map by ${\cal M}(\ell, r, t)$
called the \emph{derived map} from an admissible triple $(\ell, r,
t)$. Note that the underlying graph of the  regular map ${\cal
M}(\ell, r, t)$ is isomorphic to a coset graph ${\cal G} = {\cal
G} (\langle \ell, r, t \rangle ; \langle r, t \rangle , \langle r,
t \rangle \ell \langle r, t \rangle)$ defined as follows: $V({\cal
G}) = \{ g\langle r, t \rangle \ |\  g \in \langle \ell, r, t
\rangle \}$ and two vertices corresponding to cosets $g\langle r,
t \rangle$ and $h\langle r, t \rangle$ are adjacent if and only
if $g^{-1}h \in \langle r, t \rangle \ell \langle r, t \rangle$.
In fact, if we define $\phi: V({\cal G}) \rightarrow V(G)$ by
$\phi(g \langle r, t \rangle) = g(v)$ for any $g \in
\Gamma=\langle \ell, r, t \rangle$  then
$\phi$ is a graph isomorphism.

From now on, we consider ${\cal M}(\ell, r, t)$ as a regular
embedding of $G$ for an admissible triple $(\ell, r, t)$ for $G$.
 In \cite[Theorem 3]{GNSS}, A. Gardiner et al. showed how to
construct a regular embedding of $G$ by an admissible triple. It
looks different from our method, but they are the same in essence.

 For a combinatorial description ${\cal M}=(F;\,\lambda,\rho,\tau)$  of a regular embedding of $G$ and its corresponding admissible triple
$(\ell, r, t)$ for $G$ with a root flag $\xi \in F$,  two maps
${\cal M}(\ell, r, t)$ and ${\cal M}$ are isomorphic by an
isomorphism $\psi:{\cal M}(\ell, r, t)\rightarrow{\cal M}$ defined
by $\psi(g) = g(\xi)$ for any flag $g \in \langle \ell, r, t \rangle
= \Aut({\cal M})$.

Let  $(\ell_1,r_1,t_1)$ and $(\ell_2,r_2,t_2)$ be two admissible triples
 for $G$. If there exists a graph automorphism
$\phi$ of $G$ such that $\phi \ell_1 \phi^{-1} = \ell_2$, $\phi
r_1 \phi^{-1} = r_2$ and $\phi t_1 \phi^{-1} = t_2$ then two
derived regular maps ${\cal M}(\ell_1,r_1,t_1)$ and ${\cal
M}(\ell_2,r_2,t_2)$ are isomorphic by an isomorphism $\Phi$
defined by $\Phi(w(\ell_1, r_1, t_1)) = \phi w(\ell_1, r_1, t_1)
\phi^{-1} = w(\ell_2, r_2, t_2)$ for any word $w(\ell_1, r_1, t_1)
\in \langle \ell_1, r_1, t_1 \rangle$. Conversely, assume that two
derived regular maps ${\cal M}(\ell_1,r_1,t_1)$ and ${\cal
M}(\ell_2,r_2,t_2)$ are isomorphic by a map isomorphism $\varphi :
\langle \ell_1, r_1, t_1 \rangle \rightarrow \langle \ell_2, r_2,
t_2 \rangle$. Without loss of a generality, we can assume that
$\varphi(id)=id$. Since the underlying graphs of both
${\cal M}(\ell_1,r_1,t_1)$ and ${\cal M}(\ell_2,r_2,t_2)$ are the
same graph $G$, $\varphi$  should be  a graph automorphism of
$G$. Furthermore, it holds that $\varphi \ell_1 \varphi^{-1} =
\ell_2$, $\varphi r_1 \varphi^{-1} = r_2$ and $\varphi t_1
\varphi^{-1} = t_2$. Therefore, we get the following proposition.

 \begin{prop} \label{Thm=map-triple}
 Let $G$ be a  graph. Then,
 every regular embedding ${\cal M}$ of  $G$
into an orientable or nonorientable surface
 is
isomorphic to a derived regular map ${\cal M}(\ell ,r ,t)$  from
an admissible triple $(\ell,r,t)$ for  $G$ and its isomorphism
class corresponds to the conjugacy  class of the triple
$(\ell,r,t)$ in $\Aut(G)$.
\end{prop}

It follows from Proposition \ref{Thm=map-triple} that for a given
graph $G$, the number of nonorientable regular embeddings of $G$
up to isomorphism equals to the number of orbits of admissible
triples $(\ell,r,t)$ for $G$ satisfying $\langle \ell t, r t
\rangle = \langle \ell , r, t \rangle$ under the conjugate action
by $\Aut(G)$.

\section{ Constructions of nonorientable embeddings}

The complete bipartite graph $K_{2,2}$ is just the $4$-cycle and
there is only one nonorientable regular embedding of the $4$-cycle
with the projective plane as the supporting surface.  So, from now
on, we assume that $n \ge 3$.  For a  complete bipartite graph
$K_{n,n}$, let
  $[n] = \{ 0, 1, \ldots , n-1
 \}$ and $[n]' = \{ 0', 1', \ldots , (n-1)'
 \}$ be the
vertex sets of $K_{n,n}$  as the partite sets and let $D=\{(i,j'),
(j',i) \ | \ 0 \leq i,j \leq n-1 \}$ be the arc set.
 We denote the symmetric group on $[n] = \{
0, 1, \ldots , n-1
 \}$ by $S_{[n]}$ and the stabilizer of $i$ by $\Stab(i)$ as a subgroup of $S_{[n]}$.
 We identify the integers $0,1,\ldots, n-1$ with their
residue classes modulo $n$ according to the context.

Since $\Aut(K_{n,n})\cong S_n \wr \mathbb{Z}_2$, the wreath
product, contains all permutations of vertices of each partite set
and the interchanging of two partite sets $[n]$ and $[n]'$, one
can assume that every   orientable  or nonorientable regular
embedding of $K_{n,n}$ up to isomorphism is derived from an admissible triple $(\ell
, r_\delta , t)$  of involutions for $K_{n,n}$  of the following
type:
\begin{eqnarray*}\ell  &=& (0 \ 0') (1 \ (n-1)')(2 \
(n-2)') \cdots (n-1 \ 1') \\
 r_\delta &=& \delta (0' \ 1') ((n-1)' \ 2')((n-2)' \ 3')
\cdots \displaystyle(\left\lceil\frac{n+1}{2}\right\rceil' \
 \left\lceil\frac{n}{2}\right\rceil')
  \\
  t &=&(0)(1 \ n-1)(2 \ n-2) \cdots (\left\lfloor\frac{n}{2}\right\rfloor \
  \left\lceil\frac{n}{2}\right\rceil)(0')(1' \ (n-1)')(2' \
 (n-2)')\cdots(\left\lfloor\frac{n}{2}\right\rfloor' \
 \left\lceil\frac{n}{2}\right\rceil')
 \end{eqnarray*}
for some $\delta \in \Stab(0)$.   Note that the root vertex and
the root edge of the above admissible triple are $0$ and $(0,
0')$, respectively.  In fact, the admissibility of the triple
$(\ell , r_\delta , t)$ depends on only the permutation $\delta
\in \Stab(0)$.  Clearly, $\ell t = t \ell$ and so $\langle \ell, t
\rangle \cong \mathbb{Z}_2 \times \mathbb{Z}_2$.
 Moreover, $$ r_\delta t= \delta(0)(1 \ n-1)(2 \ n-2) \cdots (\left\lfloor\frac{n}{2}\right\rfloor \
  \left\lceil\frac{n}{2}\right\rceil)(0'\ 1' \ 2' \ \cdots
(n-2)' \ (n-1)')$$ generates the cyclic subgroup which acts regularly on
 the
arcs emanating from the root vertex $0$.

\begin{lemma} \label{delta-delta-equiv}
For any two admissible triples $(\ell, r_{\delta_1}, t)$ and
$(\ell, r_{\delta_2}, t)$   for $K_{n,n}$ with $\delta_1, \delta_2
\in \Stab(0)$, if the derived regular maps ${\cal M}(\ell,
r_{\delta_1}, t)$ and ${\cal M}(\ell, r_{\delta_2}, t)$ are
isomorphic then either {\rm (i)} $\delta_1= \delta_2$ or {\rm
(ii)} $n$ is even and $\delta_2(k) = \delta_1 (k + \frac{n}{2}) +
\frac{n}{2}$ for all $k \in [n]$.
\end{lemma}
\begin{pf}
 Suppose that ${\cal M}(\ell, r_{\delta_1},
t)$ and ${\cal M}(\ell, r_{\delta_2}, t)$ are isomorphic. By
Proposition \ref{Thm=map-triple}(2), there exists a graph
automorphism $\psi \in \Aut(K_{n,n})$ such that $\psi \ell
\psi^{-1} = \ell, \ \psi r_{\delta_1} \psi^{-1} = r_{\delta_2}$
and $\psi t \psi^{-1} = t$. Note that the vertex $0$ and possibly
$n/2$  when $n$ is even are all vertices which can be fixed by
both $t$ and $r_{\delta_i}$. The commutativity  $\psi t
 = t\psi$ and  $\psi r_{\delta_1}
 = r_{\delta_2}\psi$ implies that $\psi$ permutes these vertices.
 Hence, our discussion can be divided into the following
two cases.
\medskip

Case 1) Let $\psi(0) = 0$. Since $\psi r_{\delta_1} =
r_{\delta_2}\psi$ and $\psi$ commutes with $\ell$ and $t$, $\psi$
should be the identity. It means that $\delta_1= \delta_2$.

\medskip

Case 2) Let $n$ be even and $\psi(0) = n/2$.  Since $\psi
r_{\delta_1} = r_{\delta_2}\psi$ and $\psi$ commutes with $\ell$
and $t$, one can show that $\psi(k)=k+n/2$ and $\psi(k')=(k+n/2)'$
for all $k \in [n]$.  So, for  every  $k \in [n]$, $\ \
\delta_2(k) =r_{\delta_2}(k) = \psi r_{\delta_1} \psi^{-1}(k) =
\psi r_{\delta_1}(k+\frac{n}{2})= \delta_1 (k + \frac{n}{2}) +
\frac{n}{2}$.
\end{pf}

It will be shown later that the case (ii) in
Lemma~\ref{delta-delta-equiv} is not actually happening.

  For any involution $\delta
\in \Stab(0)$, set $$\bar{\delta} = r_\delta\cdot t|_{[n]} =
\delta \cdot(0)(1 \ -1)(2 \ -2) \cdots \in S_{[n]}.$$  Then,
$\bar{\delta}$ also belongs to $\Stab(0)$ and it satisfies an equation
$\bar{\delta}^{-1}(-k) = -\bar{\delta}(k)$ for all $k \in [n]$  as
an equivalent property of the involutory of $\delta$. For an
admissible triple $(\ell, r_{\delta}, t)$ for $K_{n,n}$ with
$\delta \in \Stab(0)$, let $R_{\bar{\delta}} = r_\delta t$ and $L
= t \ell$, namely,
\begin{eqnarray*} R_{\bar{\delta}} &=& r_\delta t =
\bar{\delta}(0' \ 1' \cdots
 (n-1)')\ \hspace{2.5cm} \mbox{and} \\
 L &=& t \ell = (0 \ 0')(1 \ 1')\cdots((n-1) \ (n-1)'),
 \end{eqnarray*}
 as permutations on the  vertex set  $[n]\cup[n]'$, which are
 the {\em local rotation automorphism} at the root vertex $0$ and the
 {\em direction-reversing automorphism of the root edge $(0\ 0')$}, respectively.
  In fact, $L$ is an automorphism which interchange partite sets.
 Note that any one of $\delta, \bar{\delta}, R_{\bar{\delta}}$  determines
 completely the other two.
 Recall that the regular map ${\cal M}(\ell, r_{\delta},
t)$ is nonorientable if and only if $\left< \ell, r_{\delta}, t
\right> = \left< R_{\bar{\delta}},L \right>$. Hence, if $(\ell,
r_{\delta}, t)$ is an admissible triple for $K_{n,n}$ and ${\cal
M}(\ell, r_{\delta}, t)$ is nonorientable then $\bar{\delta}(0)=0$,
$\bar{\delta}^{-1}(-k) = -\bar{\delta}(k)$ for all $k \in [n]$,
$|\left< R_{\bar{\delta}}, L \right>|=4|E(K_{n,n})|= 4n^2$ and $t
\in \left< R_{\bar{\delta}},L \right>$. So, in order to construct
all nonorientable regular embeddings of $K_{n,n}$, we need to
examine $\bar{\delta}$ satisfying the aforementioned conditions. Let
\begin{eqnarray*} \mathcal{M}^{non}_n &=&\{ \bar{\delta} \in S_{[n]} \
| \  \bar{\delta}(0)=0, \bar{\delta}^{-1}(-k) = -\bar{\delta}(k) \ \ \ \mbox{for all
$k \in [n]$}
 \\ & & \hspace{3cm} \ |\left< R_{\bar{\delta}}, L \right>|= 4n^2 \ \ \mbox{and $\left< R_{\bar{\delta}},L \right>$ contains $t$} \}. \end{eqnarray*}

 We shall show that there is a one-to-one
correspondence between the nonorientable regular embeddings of
$K_{n,n}$ for $n \geq 3$ up to isomorphism and the elements in
$\mathcal{M}^{non}_n$. From now on, we shall deal with
$\bar{\delta}$ instead of $\delta$ to construct and to classify the
nonorientable regular embeddings of $K_{n,n}$.

\begin{lemma} \label{alpha-sigma-equiv}
 For every  involution $\delta \in S_{[n]}$ with $\delta(0)=0$, the
following statements are equivalent.
\begin{enumerate}
\item[$(1)$]  The triple $(\ell, r_{\delta}, t)$ is admissible and
the derived regular map ${\cal M}(\ell, r_{\delta}, t)$ is
nonorientable.
 \item[$(2)$] $\bar{\delta}\in
 \mathcal{M}^{non}_n$, where $\bar{\delta}=r_\delta\cdot
 t|_{[n]}.$
\end{enumerate}
\end{lemma}

\begin{pf}  For $\bar{\delta}=r_\delta\cdot
 t|_{[n]},$ we know already  $\bar{\delta}(0)=0$ and $\bar{\delta}^{-1}(-k) =
-\bar{\delta}(k)$ for all $k \in [n]$.

$(1)\Rightarrow (2)$\ Let $(\ell, r_{\delta}, t)$ be admissible
and let the map ${\cal M}(\ell, r_{\delta}, t)$ be nonorientable.
Then, $\left< \ell, r_{\delta}, t \right> = \left<
R_{\bar{\delta}},L \right>$ and
 $|\left< \ell, r_{\delta}, t \right>|=|\left< R_{\bar{\delta}}, L
\right>|=4|E(K_{n,n})|= 4n^2$. So, $\bar{\delta} \in
\mathcal{M}^{non}_n$.

$(2)\Rightarrow (1)$\ Let $\bar{\delta} \in \mathcal{M}^{non}_n$.
Since $t \in \left< R_{\bar{\delta}},L \right>$, it holds that
$R_{\bar{\delta}} t = r_\delta \in \left< R_{\bar{\delta}},L
\right>$ and
 $t L = \ell \in \left< R_{\bar{\delta}},L \right>$. Hence,
 $\left< \ell, r_{\delta}, t \right> = \left< R_{\bar{\delta}},L \right>$.
 For any $i,j \in [n]$, we have
 $$R_{\bar{\delta}}^{i} L  R_{\bar{\delta}}^{j} t
 (0,0') = R_{\bar{\delta}}^{i} L  R_{\bar{\delta}}^{j}
 (0,0') = R_{\bar{\delta}}^{i} L (0,j') = R_{\bar{\delta}}^{i}(0', j)
 = (i',\bar{\delta}^i(j))$$
and by taking $L$ on both sides, we have
 $$L R_{\bar{\delta}}^{i} L  R_{\bar{\delta}}^{j} t
 (0,0') = L R_{\bar{\delta}}^{i} L  R_{\bar{\delta}}^{j}
 (0,0') = (i,\bar{\delta}^i(j)').$$
This shows that the arc $(0,0')$ can be mapped to  any other arc  by the action
of  the group $\left< R_{\bar{\delta}},L \right>$. It means that
$\left< \ell, r_{\delta}, t \right> = \left< R_{\bar{\delta}},L
\right>$ acts transitively on both the arc set  $D(K_{n,n})$ and the
vertex set
 $V(K_{n,n})$. For $0 \in V(K_{n,n})$, $\left< R_{\bar{\delta}},t
\right> \leq \left< R_{\bar{\delta}},L \right>_0$. Since $|\left<
R_{\bar{\delta}},t \right>| = | \left< R_{\bar{\delta}},L
\right>_0|=2n$, one can see that $\left< R_{\bar{\delta}},t
\right> = \left< R_{\bar{\delta}},L \right>_0 \simeq D_{n}$ of
order $2n$ and the subgroup $\left< R_{\bar{\delta}} \right>$ acts
regularly on the arcs emanating from $0$. For the edge
$e=\{0,0'\}$, one can easily check that the stabilizer $\left<
R_{\bar{\delta}},L \right>_e$ is equal to $\left< L,t \right>$
which is isomorphic to $\mathbb{Z}_2 \times \mathbb{Z}_2$. So,
$(\ell, r_{\delta}, t)$ is an admissible triple for $K_{n,n}$.
Since $\left< \ell, r_{\delta}, t \right> = \left<
R_{\bar{\delta}},L \right>$, the derived regular map ${\cal
M}(\ell, r_{\delta}, t)$ is nonorientable.
\end{pf}

To determine $\bar{\delta}$ satisfying $|\left< R_{\bar{\delta}},L
\right>|=4n^2$, we need to examine the words of $R_{\bar{\delta}}$
and $L$ in the group $\left< R_{\bar{\delta}},L \right>$.

\begin{lemma} \label{disjoint-union}
 Let  $\delta \in S_{[n]}$  be an involution and
$\delta(0)=0$, or equivalently $\bar{\delta}^{-1}(-k) =
-\bar{\delta}(k)$ for all $k \in [n]$ and $\bar{\delta}(0)=0$.
Then the following statements are equivalent.
\begin{enumerate}
\item[$(1)$]  $\bar{\delta} \in \mathcal{M}^{non}_n$.
 \item[$(2)$]        The subgroup $\left<R_{\bar{\delta}}, L \right> $ of $S_{[n]\cup[n]'}$ is a disjoint union of the four sets
\begin{eqnarray*}
  B &:=& \{R_{\bar{\delta}}^{i} L
R_{\bar{\delta}}^{j} \ | \ i,j \in [n] \}, \\
  LB &:=& \{ L R_{\bar{\delta}}^{i} L
R_{\bar{\delta}}^{j} \ | \ i,j \in [n] \}, \\
  Bt &:=& \{  R_{\bar{\delta}}^{i} L R_{\bar{\delta}}^{j} t \ | \
i,j \in [n] \}, \\
  LBt &:=& \{ L R_{\bar{\delta}}^{i} L R_{\bar{\delta}}^{j} t \ | \
i,j \in [n] \}.
\end{eqnarray*}
 \item[$(3)$]  For each $i \in [n]$,  there exist
 $a(i), b(i) \in [n]$ such that  either \vspace{-2mm}
 \[ \bar{\delta}(k + i) =
 \bar{\delta}^{b(i)}(k) + a(i) \ \ \mbox{and} \ \ \bar{\delta}^{i}(k)+1 =
 \bar{\delta}^{a(i)}(k + b(i)) \ \ \mbox{for all $k \in [n],$} \hspace{1cm} \tag{$\ast_1$} \vspace{-2mm}\]
 or \vspace{-2mm} \[ \bar{\delta}(k + i) =
 \bar{\delta}^{b(i)}(-k) + a(i) \ \ \mbox{and} \ \ \bar{\delta}^{i}(k)+1 =
 \bar{\delta}^{a(i)}(-k + b(i))\ \ \mbox{for all $k \in [n]$}. \hfill \tag{$\ast_2$}\]
  In addition, the latter case Eq.$(\ast_2)$ holds for at
least one $i \in [n]$.
\end{enumerate}
\end{lemma}

 \begin{pf} $(1) \Rightarrow (2)$ \
 For any $\bar{\delta} \in \mathcal{M}^{non}_n$ and for any $i,j \in [n]$, we have
 $$R_{\bar{\delta}}^{i} L  R_{\bar{\delta}}^{j} t
 (0,0') = R_{\bar{\delta}}^{i} L  R_{\bar{\delta}}^{j}
 (0,0') = R_{\bar{\delta}}^{i} L (0,j') = R_{\bar{\delta}}^{i}(0', j)
 = (i',\bar{\delta}^i(j))$$
and by taking $L$ on both sides
 $$L R_{\bar{\delta}}^{i} L  R_{\bar{\delta}}^{j} t
 (0,0') = L R_{\bar{\delta}}^{i} L  R_{\bar{\delta}}^{j}
 (0,0') = (i,\bar{\delta}^i(j)').$$
 By comparing images of the arc $(0,0')$, one can see that $$ \left( B \cup Bt \right) \cap \left( LB \cup LBt \right) = \emptyset$$ and for
any $(i,j) \neq (k, \ell)$, it holds that $ \{
R_{\bar{\delta}}^{i} L R_{\bar{\delta}}^{j}, R_{\bar{\delta}}^{i}
L R_{\bar{\delta}}^{j}t \} \cap \{ R_{\bar{\delta}}^{k} L
R_{\bar{\delta}}^{\ell}, R_{\bar{\delta}}^{k} L
R_{\bar{\delta}}^{\ell}t \} = \emptyset$ and $\{
LR_{\bar{\delta}}^{i} L R_{\bar{\delta}}^{j},
LR_{\bar{\delta}}^{i} L R_{\bar{\delta}}^{j}t  \} \cap \{
LR_{\bar{\delta}}^{k} L R_{\bar{\delta}}^{\ell},
LR_{\bar{\delta}}^{k} L R_{\bar{\delta}}^{\ell}t  \} = \emptyset$.
Now, it suffices to show that $R_{\bar{\delta}}^{i} L
R_{\bar{\delta}}^{j}\neq R_{\bar{\delta}}^{i} L
R_{\bar{\delta}}^{j}t$ and $LR_{\bar{\delta}}^{i} L
R_{\bar{\delta}}^{j}\neq LR_{\bar{\delta}}^{i} L
R_{\bar{\delta}}^{j}t$ for all $(i,j) \in [n]\times [n]$ for
disjointness of $B \cap Bt = \emptyset$ and $LB \cap LBt =
\emptyset$. In fact, for any $(i,j) \in [n]\times [n]$,
$$R_{\bar{\delta}}^{i} L R_{\bar{\delta}}^{j}t
 (0,1') = R_{\bar{\delta}}^{i} L  R_{\bar{\delta}}^{j}(0,-1')$$
and
$$LR_{\bar{\delta}}^{i} L  R_{\bar{\delta}}^{j}t
 (0,1') = LR_{\bar{\delta}}^{i} L  R_{\bar{\delta}}^{j}
 (0,-1').$$
 These implies that $ R_{\bar{\delta}}^{i} L
R_{\bar{\delta}}^{j} \neq  R_{\bar{\delta}}^{i} L
R_{\bar{\delta}}^{j}t$ and $LR_{\bar{\delta}}^{i} L
R_{\bar{\delta}}^{j} \neq LR_{\bar{\delta}}^{i} L
R_{\bar{\delta}}^{j}t$. Hence, the four sets $B, LB, Bt$ and $LBt$
are mutually disjoint and the cardinality of their union is
$4n^2$, which equals $|\left<R_{\bar{\delta}}, L \right> |$.

$(2) \Rightarrow (3)$ \ Since the group $\left<R_{\bar{\delta}}, L
\right>$ is the union of the four sets, for each $i \in [n]$, there
exist
 $a(i), b(i) \in [n]$ such that $R_{\bar{\delta}} L
R_{\bar{\delta}}^{i} L = L R_{\bar{\delta}}^{a(i)} L
R_{\bar{\delta}}^{b(i)}$ or $R_{\bar{\delta}} L
R_{\bar{\delta}}^{i} L = L R_{\bar{\delta}}^{a(i)} L
R_{\bar{\delta}}^{b(i)}t$. By comparing their values of $k$ and
$k'$, we have   $$\bar{\delta}(k + i) =
 \bar{\delta}^{b(i)}(k) + a(i) \ \ \mbox{and}\ \ \bar{\delta}^{i}(k)+1 =
 \bar{\delta}^{a(i)}(k + b(i)) \ \ \mbox{for all $k \in [n]$} $$ or $$\bar{\delta}(k + i) =
 \bar{\delta}^{b(i)}(-k) + a(i) \ \ \mbox{and}\ \ \bar{\delta}^{i}(k)+1 =
 \bar{\delta}^{a(i)}(-k + b(i)) \ \ \mbox{for all $k \in [n]$}.$$
  That is, either Eq.$(\ast_1)$ or  Eq.$(\ast_2)$ holds.
 Suppose that Eq.$(\ast_1)$ holds for
all $i \in [n]$,  namely, $R_{\bar{\delta}} L R_{\bar{\delta}}^{i} L
= L R_{\bar{\delta}}^{a(i)} L R_{\bar{\delta}}^{b(i)}.$ Then  one
can easily check that
$$\left<R_{\bar{\delta}}, L \right> = \{ L R_{\bar{\delta}}^{i} L
R_{\bar{\delta}}^{j} \ | \ i,j \in [n] \} \ \cup \ \{
R_{\bar{\delta}}^{i} L R_{\bar{\delta}}^{j} \ | \ i,j \in [n] \} \
\ \mbox{\rm (disjoint union),}$$ which contradicts the assumption.
Hence, there exists at least one $i \in [n]$ for which
Eq.$(\ast_2)$ holds.

 $(3) \Rightarrow (1)$ \ Note that Eq.$(\ast_1)$ is
 nothing but  the equality
  $R_{\bar{\delta}} L
R_{\bar{\delta}}^{i} L = L R_{\bar{\delta}}^{a(i)} L
R_{\bar{\delta}}^{b(i)}$ and   Eq.$(\ast_2)$ is
 equivalent to the equality
  $R_{\bar{\delta}} L
R_{\bar{\delta}}^{i} L = L R_{\bar{\delta}}^{a(i)} L
R_{\bar{\delta}}^{b(i)}t$. These two equalities imply that
$\left<R_{\bar{\delta}}, L \right>$ contains $t$ and
$\left<R_{\bar{\delta}}, L \right>$ is the same as  the union of
the four sets in $(2)$. The disjointness of the union can be shown
in a similar way as in $(1) \Rightarrow (2)$. It means that
$|\left<R_{\bar{\delta}}, L \right> | = 4n^2$. So, $\bar{\delta}
\in \mathcal{M}^{non}_n$.
\end{pf}

In fact, two numbers $a(i)$ and $b(i)$ in
Lemma~\ref{disjoint-union}$(3)$ are completely determined by the
values of $\bar{\delta}$ in the following sense.

\begin{lemma} \label{a(i)b(i)}
 Let $\bar{\delta} \in \mathcal{M}^{non}_n$ and let $a(i)$ and
 $b(i)$ be the numbers given in Eqs.$(\ast_1)$ and $(\ast_2)$. If they
 satisfy Eq.$(\ast_1)$ then \vspace{-2mm}
 \[ a(i)=\bar{\delta}(i) \ \ \mbox{and} \ \
 b(i)=\bar{\delta}^{i}(1)=\bar{\delta}^{-\bar{\delta}(i)}(1), \vspace{-2mm} \]
 and if they
 satisfy Eq.$(\ast_2)$ then \vspace{-2mm}
 \[ a(i)=\bar{\delta}(i) \ \ \mbox{and} \ \
 b(i)=-\bar{\delta}^{i}(1)=\bar{\delta}^{-\bar{\delta}(i)}(1). \vspace{-2mm} \]
\end{lemma}
\begin{pf}\ Let $a(i)$ and
 $b(i)$  satisfy Eq.$(\ast_1)$.  By taking $k=0$ in the equation $\bar{\delta}(k + i) =
 \bar{\delta}^{b(i)}(k) + a(i)$, we have
 $a(i)=\bar{\delta}(i)$. And, by taking $k=0$ and $k=-b(i)$ in
 the equation $\bar{\delta}^{i}(k)+1 = \bar{\delta}^{a(i)}(k + b(i))
 $, we have
 $b(i)=\bar{\delta}^{-a(i)}(1)=\bar{\delta}^{-\bar{\delta}(i)}(1)$
 and $b(i)=-\bar{\delta}^{-i}(-1)$. Since $\bar{\delta}^{-1}(-k)
 = -\bar{\delta}(k)$ for all $k \in [n]$,
 \begin{eqnarray*} b(i)&=&-\bar{\delta}^{-i}(-1)=-\bar{\delta}^{-i+1}(\bar{\delta}^{-1}(-1))=-\bar{\delta}^{-i+1}(-\bar{\delta}(1))
=-\bar{\delta}^{-i+2}(-\bar{\delta}^2(1))= \cdots \\
&=&-\bar{\delta}^{-1}(-\bar{\delta}^{i-1}(1))=\bar{\delta}^{i}(1).\end{eqnarray*}

Let $a(i)$ and
 $b(i)$  satisfy Eq.$(\ast_2)$. By taking $k=0$ in the equation $\bar{\delta}(k
+ i) =
 \bar{\delta}^{b(i)}(-k) + a(i)$, one has
 $a(i)=\bar{\delta}(i)$. And, by taking $k=0$ and $k=b(i)$ in
 the equation $\bar{\delta}^{i}(k)+1 = \bar{\delta}^{a(i)}(-k + b(i))
 $, one can see that
 $b(i)=\bar{\delta}^{-a(i)}(1)=\bar{\delta}^{-\bar{\delta}(i)}(1)$
 and $b(i)=\bar{\delta}^{-i}(-1)$. Since $\bar{\delta}^{-1}(-k)
 = -\bar{\delta}(k)$ for all $k \in [n]$, $b(i)=\bar{\delta}^{-i}(-1) = -\bar{\delta}^{i}(1)$.
\end{pf}

 Now, let us consider the even numbers $n$ as the
first case to construct nonorientable regular embeddings of
$K_{n,n}$. As a candidate of $\bar{\delta} \in
\mathcal{M}^{non}_n$, we define  a permutation $\bar{\delta}_{n,x}
\in S_{[n]}$  by
\[ \ \bar{\delta}_{n,x}=(0)(2 \ \ -2)(4 \ \ -4)\cdots(1 \ \ 1+x \ \ 1+2x \ \ 1+3x \ \
\cdots) \ \]  for any positive  even integer $x$ such that the
greatest common divisor of $n$ and $x$ is 2. In fact, the $\bar{\delta}_{n,x}$'s are all possible permutations in $\mathcal{M}^{non}_n$ which do not have a reduction which will be defined in the next section.

Suppose that $\bar{\delta}_{n,x}$ belongs to $\mathcal{M}^{non}_n$ for some even $x$. Then, for each $i \in [n]$ there
exist  $a(i), b(i) \in [n]$ satisfying Eq.$(\ast_1)$ or  Eq.$(\ast_2)$  by Lemma~\ref{disjoint-union}. For every even $i \in [n]$, it holds that
\[ b(i)= \bar{\delta}_{n,x}^{-\bar{\delta}_{n,x}(i)}(1) = \bar{\delta}_{n,x}^{i}(1), \]
which implies that $a(i),b(i)$ satisfy Eq.$(\ast_1)$ by Lemma~\ref{a(i)b(i)}. So, there exists an odd integer $i \in [n]$ such that $a(i), b(i) \in [n]$ satisfying Eq.$(\ast_2)$ by Lemma~\ref{disjoint-union}. For such odd $i \in [n]$,
\[ b(i) =- \bar{\delta}_{n,x}^{i}(1) = -(1+ix) \ \ \mbox{and} \]
  \[ b(i)=\bar{\delta}_{n,x}^{-\bar{\delta}_{n,x}(i)}(1) = \bar{\delta}_{n,x}^{-(i+x)}(1) = 1-(i+x)x=-(1+ix)+(2-x^2) \]
by Lemma~\ref{a(i)b(i)}. It implies that $x^2 \equiv 2 (\mod n)$. Therefore, the condition $x^2 \equiv 2 (\mod n)$ is necessary for $\bar{\delta}_{n,x}$ to belong to $\mathcal{M}^{non}_n$. The following two lemmas show that the condition is also sufficient.

\begin{lemma} \label{check-admissible}
 Let $n$ and $x$ be even integers such that $n > x > 3$ and
$x^2 \equiv 2 (\mod n)$.
 \begin{enumerate} \item[$(1)$] For an even integer $2i
\in [n],$ if we
define $a(2i) = -2i$ and $b(2i) = 2ix+1$  then
$a(2i)$ and $b(2i)$ satisfy Eq.$(\ast_1)$ with
$\bar{\delta}_{n,x}$.
 \item[$(2)$] For an odd integer $2i+1 \in [n],$ if we define $a(2i+1) = 2i+x+1$
and $b(2i+1) = -2ix-x-1$  then
$a(2i+1)$ and $b(2i+1)$ satisfy Eq.$(\ast_2)$ with
$\bar{\delta}_{n,x}$.
\end{enumerate}
\end{lemma}
\begin{pf} Since a proof is similar,  we prove only $(1).$ \ For an even integer $2i,$ let
$a(2i) = -2i$ and $b(2i) = 2ix+1.$ Then for any even integer $2k
\in [n]$, we have
 \begin{eqnarray*} \bar{\delta}_{n,x}(2k+2i) &=& -2k-2i \\
\bar{\delta}_{n,x}^{b(2i)}(2k)+a(2i) &=&
\bar{\delta}_{n,x}^{2ix+1}(2k)-2i
=-2k -2i  \\
 \bar{\delta}_{n,x}^{2i}(2k)+1&=& 2k+1 \ \ \mbox{and} \\
\bar{\delta}_{n,x}^{a(2i)}(2k+b(2i)) &=&
\bar{\delta}_{n,x}^{-2i}(2k+2ix+1)=2k+2ix+1 -2ix =2k+1.
 \end{eqnarray*}
 For any odd integer $2k+1 \in [n]$, we have
\begin{eqnarray*} \bar{\delta}_{n,x}(2k+1+2i) &=& 2k+2i+x+1 \\
\bar{\delta}_{n,x}^{b(2i)}(2k+1)+a(2i) &=&
\bar{\delta}_{n,x}^{2ix+1}(2k+1)-2i
=2k+1 +(2ix+1)x -2i \\
&=& 2k+1+4i+x-2i = 2k+2i+x+1 \  \mbox{because} \   x^2 \equiv 2 (\mod n)\\
 \bar{\delta}_{n,x}^{2i}(2k+1)+1&=& 2k+1+2ix+1 =2k+2ix+2\ \ \mbox{and} \\
\bar{\delta}_{n,x}^{a(2i)}(2k+1+b(2i)) &=&
\bar{\delta}_{n,x}^{-2i}(2k+2ix+2)=2k+2ix+2.
 \end{eqnarray*}
 Hence, $a(2i)$ and $b(2i)$ satisfy Eq.$(\ast_1)$ with
the
 permutation
$\bar{\delta}_{n,x}$.
\end{pf}

\begin{lemma}  \label{construction}
For any two positive even integers $n$ and $x$ such that  $n > x > 3$ and $x^2 \equiv 2 (\mod n)$,
$\bar{\delta}_{n,x}$ belongs to $\mathcal{M}^{non}_n$.
\end{lemma}

\begin{pf}
Note that $\bar{\delta}_{n,x}(0)=0$. For any even $2k \in [n]$,
$\bar{\delta}_{n,x}^{-1}(-2k) = 2k = -\bar{\delta}_{n,x}(2k)$. For
any odd $2k+1 \in [n]$, $\bar{\delta}_{n,x}^{-1}(-2k-1) = -2k-1 -
x$ and $ - \bar{\delta}_{n,x}(2k+1)=-(2k+1+x) = -2k-1 - x$. So,
for any $k \in [n]$, it holds that $\bar{\delta}_{n,x}^{-1}(-k) =
- \bar{\delta}_{n,x}(k)$. By Lemmas~\ref{disjoint-union}(3) and
 \ref{check-admissible}, $\bar{\delta}_{n,x}$ belongs to
$\mathcal{M}^{non}_n$.
\end{pf}

By Lemma~\ref{construction}, as a subset of $\mathcal{M}^{non}_n$ let us consider the following set:
 \[ \mathcal{N}^{non}_n = \left\{ \begin{array}{cl} \{\bar{\delta}_{n,x} \ | \ n
> x > 3 \ \ \mbox{and} \ \ x^2 \equiv 2 (\mod n)  \ \} & \mbox{ if $n$ is even} \\
\emptyset
& \mbox{ if $n$ is odd}
\end{array}. \right. \]

Note that for any even integers $n, x$ such that $n \equiv 0 (\mod 4)$ and $n > x > 1$, $x^2$ is a multiple of 4 and hence there is no such $x$ satisfying $x^2 \equiv 2 (\mod n)$. So,  for every $n \equiv 0 (\mod 4)$,
$\mathcal{N}^{non}_n = \emptyset$. The smallest integer $n$ such
that $\mathcal{N}^{non}_n \neq \emptyset$ is $14$.  This will show
that $K_{14,14}$ is the smallest complete bipartite graph which
can be regularly embedded into a nonorientable surface except
$K_{2,2}$.

In the last two sections, it will be shown that
 $\mathcal{M}^{non}_n = \mathcal{N}^{non}_n$ for every $n$,
 which means that $\mathcal{M}^{non}_n =\emptyset$ if $n$ is  odd or
  $n \equiv 0 (\mod 4)$.
 In the remaining of
this section, we shall show that for any two different
 $\bar{\delta}_{n,x_1}, \bar{\delta}_{n,x_2} \in
\mathcal{N}^{non}_n$, their derived regular embeddings of
$K_{n,n}$ are not isomorphic. And, for a given $n \equiv 2 (\mod
4)$, we will estimate the cardinality $|\mathcal{N}^{non}_n|$,
that is, the number of solutions of $x^2 = 2$ in $\mathbb{Z}_n$.

\begin{lemma} \label{admissible-pair-isomor}
For any two $\bar{\delta}_{n,x_1}, \bar{\delta}_{n,x_2} \in
\mathcal{N}^{non}_n$ with $n > 3$, let
$\delta_i=\bar{\delta}_{n,x_i}\cdot (0)(1 \ -1)(2 \ -2)\cdots$ for
$i=1,2$. Then, two derived regular maps ${\cal M}(\ell,
r_{\delta_1}, t)$ and ${\cal M}(\ell, r_{\delta_2}, t)$ are
isomorphic if and only if $x_1 = x_2$.
\end{lemma}
\begin{pf}
Since the sufficiency is clear, we prove only the necessity.
Recall that if $n \equiv 0 (\mod 4)$ then $\mathcal{N}^{non}_n =
\emptyset$. So, let $n \equiv 2 (\mod 4)$.

Let two  regular maps ${\cal M}(\ell, r_{\delta_1}, t)$ and ${\cal
M}(\ell, r_{\delta_2}, t)$ be isomorphic. By
Lemma~\ref{delta-delta-equiv}, $\delta_1 = \delta_2$ or
$\delta_2(k) = \delta_1(k + \frac{n}{2})+ \frac{n}{2}$ for any $k
\in [n]$. If $\delta_1 = \delta_2$ then $x_1 =x_2$. Suppose that
$\delta_2(k) =\delta_1(k + \frac{n}{2})+ \frac{n}{2}$ for any $k
\in [n]$. By taking $k=0$ in the equation $\delta_2(k) =\delta_1(k
+ \frac{n}{2})+ \frac{n}{2}$, one can get
\[ \ 0=\delta_2(0) = \delta_1(\frac{n}{2})+\frac{n}{2} =
\bar{\delta}_{n,x_1}(\frac{n}{2})+\frac{n}{2} =
\frac{n}{2}+x_1+\frac{n}{2}=x_1. \] Since $x_1^2 \equiv 2 (\mod n)$,
 this is impossible.
\end{pf}

The following two lemmas are well-known in number theory. So, we
state them without a proof. (See p.112 and p.77 of the book
\cite{AG}.)

\begin{lemma} \label{Gauss-lemma} {\bf (Gauss' lemma)} Let $p$ be
an odd prime and let $a$ be an integer such that $p \nmid a$.
Consider a sequence of integers $a,2a,3a,\ldots,(\frac{p-1}{2})a$.
Replace each integer in the sequence by the one congruent to it
modulo $p$ which lies between $-\frac{p-1}{2}$ and
$\frac{p-1}{2}$. Let $\nu$ be the number of negative integers in
the resulting sequence. Then, $x^2 \equiv a (\mod p)$ has a
solution if and only if $\nu$ is even.
\end{lemma}

\begin{corollary} \label{is2quadreci?}
For any odd prime $p$, $x^2 \equiv 2 (\mod p)$ has a solution if
and only if $p \equiv \pm 1 (\mod 8)$.
\end{corollary}

\begin{lemma} \label{p-p^n}  Let $p$ be
an odd prime and let $a$ be an integer such that $p \nmid a$. Then,
for any positive integer $m$, $x^2 \equiv a (\mod p)$ has a solution
 in $\mathbb{Z}_p$ if and
only if $x^2 \equiv a (\mod p^m)$ has a solution  in
$\mathbb{Z}_{p^m}$. Moreover, they have the same number of
solutions, which is 0 or 2.
\end{lemma}

 Since  $\mathcal{N}^{non}_n = \emptyset$ for $n \equiv 0 (\mod 4)$, we need to
 estimate $|\mathcal{N}^{non}_n|$ for only
 $n \equiv 2 (\mod 4)$.

\begin{lemma} \label{card-Mnmp}
 For an $n =
2p_1^{a_1}p_2^{a_2}\cdots p_k^{a_k}$  (a prime decomposition), the
number $|\mathcal{N}^{non}_n|$ of solutions of $x^2 =2$ in
$\mathbb{Z}_n$ is $2^k$ if $p_i \equiv \pm 1 (\mod 8)$ for all
$i=1,2,\ldots,k$; 0 otherwise.
\end{lemma}
\begin{pf}
For an $x \in \mathbb{Z}_n$, $x^2 \equiv 2 (\mod n)$ if and only
if $x$ is even and $x^2 \equiv 2 (\mod p_i^{a_i})$ for all
$i=1,2,\ldots,k$. Hence, by Corollary \ref{is2quadreci?} and Lemma
\ref{p-p^n}, if  $p_i \equiv \pm 3 (\mod 8)$
 for some $i \geq 1$, the
cardinality $|\mathcal{N}^{non}_n|$ is zero. If $p_i \equiv \pm 1
(\mod 8)$ for all $i=1,2,\ldots,k,$ $|\mathcal{N}^{non}_n|=2^k$ by
Corollary \ref{is2quadreci?}, Lemma \ref{p-p^n} and the Chinese
remainder theorem.
\end{pf}

\section{Reduction}

In this section, we show that if there would exist a $\bar{\delta}
\in \mathcal{M}^{non}_n - \mathcal{N}^{non}_n$, then $d=| \langle
\bar{\delta} \rangle | < n$ and there is an induced element
$\bar{\delta}_{(1)}$ in $\mathcal{M}^{non}_d$, called the
\emph{reduction} of $\bar{\delta}$. If such $\bar{\delta}_{(1)}$
is also contained in $\mathcal{M}^{non}_d - \mathcal{N}^{non}_d$
then one can choose the next reduction $\bar{\delta}_{(2)}$ of
$\bar{\delta}_{(1)}$. By continuing such reduction, one can have a
nonnegative integer $j$ such that $\bar{\delta}_{(j)} \in
\mathcal{M}^{non}_{d_j} - \mathcal{N}^{non}_{d_j}$  but its
reduction $\bar{\delta}_{(j+1)}$ is the identity or belongs to $
\mathcal{N}^{non}_{d_{j+1}}$. In the next section, we prove
$\mathcal{M}^{non}_n = \mathcal{N}^{non}_n$ for any $n$ by showing
that such a $\bar{\delta}_{(j)}$ does not exist.

For an even integer $n$, let  us write  $\bar{n}=n/2$ for
notational convenience.
 The following lemma  is related to the order
of $\bar{\delta} \in \mathcal{M}^{non}_n$.

\begin{lemma} \label{Lem=order-sigma}
 Suppose that $\bar{\delta} \in \mathcal{M}^{non}_n -\mathcal{N}^{non}_n$ exists.
 Then, the order of the cyclic group $\left<\bar{\delta}\right>$
 equals the size of the orbit of
 $1$ under $\langle \bar{\delta}\rangle$, namely, $|\{\bar{\delta}^{i}(1)\ | \ i \in [n]
\}|$. Furthermore, it is a divisor of $n$ but not equal to $n$.
\end{lemma}

\begin{pf}
 For any $k \in [n]$, let $O(k) = \{\bar{\delta}^{i}(k) \ | \ i \in [n] \}$ be the orbit of  $k$ under $\langle \bar{\delta}\rangle$.
 Let  $|O(1)| =
 d$. Then, $d$ divides $n$ and  $d < n$ because $0 \notin O(1)$. Moreover, we have \vspace{-.2cm}
$$\bar{\delta}^{d}(1) = 1~~\mbox{and} ~~(L R_{\bar{\delta}}
L)^{-1} R_{\bar{\delta}}^{d} (L R_{\bar{\delta}} L) (0) =
0.\vspace{-.2cm}$$ Hence,  the conjugate $(L R_{\bar{\delta}}
L)^{-1} R_{\bar{\delta}}^{d}( L R_{\bar{\delta}} L)$ of
$R_{\bar{\delta}}^d$ belongs to the vertex stabilizer $
\left<R_{\bar{\delta}}, L \right>_{0} = \left<R_{\bar{\delta}}, t
\right>$ which is isomorphic to a dihedral group $D_n$ of order
$2n$.

Assume that  $(L R_{\bar{\delta}} L)^{-1} R_{\bar{\delta}}^{d}( L
R_{\bar{\delta}} L) = R_{\bar{\delta}}^{m}$ for some $m \in [n]$.
 Because $R_{\bar{\delta}}^{m}$ and
$R_{\bar{\delta}}^{d}$ are conjugate in $\left<R_{\bar{\delta}},
L\right>$, we have
$\left<R_{\bar{\delta}}^{m}\right>=\left<R_{\bar{\delta}}^{d}\right>$
as subgroups of the cyclic group $\left<R_{\bar{\delta}}\right>$.
Since $d$ is a divisor of $n$, there exists $ \ell \in
[\frac{n}{d}]$ such that $m = \ell d$ and $( \ell , \frac{n}{d}) =
1$.
 Suppose  that $|\left<\bar{\delta}\right>| \neq
d$. Then, there exists $k \in [n]$ such that $\bar{\delta}^{d} (k)
\neq k$. Let $q$ be the largest  such $k$. Then, $\bar{\delta}^{
\ell  d} (q) \neq q$. On the other hand,
$$\bar{\delta}^{ \ell  d} (q)=R_{\bar{\delta}}^{ \ell d}(q) =R_{\bar{\delta}}^{m}(q)
=(L R_{\bar{\delta}} L)^{-1} R_{\bar{\delta}}^{d} (L
R_{\bar{\delta}} L) (q) = (L R_{\bar{\delta}} L)^{-1}
R_{\bar{\delta}}^{d} (q+1) = q,$$ which contradicts
$\bar{\delta}^{ \ell d} (q) \neq q$. Therefore,
$|\left<\bar{\delta}\right>| = |O(1)| = d$, a divisor of $n$ but
$d \neq n.$

Next, suppose that  $(L R_{\bar{\delta}} L)^{-1}
R_{\bar{\delta}}^{d}( L R_{\bar{\delta}} L) =
R_{\bar{\delta}}^{m}t$ for some $m \in [n]$. Then, since the order
of $R_{\bar{\delta}}^{m}t$ is 2 and $d < n$, $n$ is even and
$d=\bar{n}=n/2$. If $|O(k)|$ divides $\bar{n}$ for all $k \in [n]$
then $\left<\bar{\delta}\right>$ is a cyclic group of order
$\bar{n}$ and the result follows. So, we may assume that there is
some $i \in [n]$ such that $|O(i)|$  doesn't divide $\bar{n}$. By
comparing two values $(L R_{\bar{\delta}} L)^{-1}
R_{\bar{\delta}}^{\bar{n}}( L R_{\bar{\delta}} L)(k)$ and
$R_{\bar{\delta}}^{m}t(k)$, we have
$$\bar{\delta}^{\bar{n}}(k+1)-1 = \bar{\delta}^m(-k) \ \ \mbox{or equivalently} \ \ \bar{\delta}^{\bar{n}}(k+1) = \bar{\delta}^m(-k)+1 $$ for all
$k \in [n]$. Note that if $\bar{\delta}^{\bar{n}}(k+1) = k+1$ for
some $k \in [n]$,
 then $\bar{\delta}^m(-k) = k$.
 Let $\Fix(\bar{\delta})=\{ k \in [n] \ | \ \bar{\delta}(k)=k \}$ and $f(\bar{\delta}) =
 |\Fix(\bar{\delta})|$. Then, $f(\bar{\delta})
 \geq 1$ by the fact $0 \in \Fix(\bar{\delta})$.
Since $\bar{\delta}^{\bar{n}}(j) = j$ for any $j \in O(1) \cup
\Fix(\bar{\delta})($disjoint union$)$, there are at least $\bar{n}
+ f(\bar{\delta}) \geq \bar{n} +1$ elements $k \in [n]$ satisfying
$\bar{\delta}^m(-k) = k$ or equivalently, there exist at most $\bar{n} -1$ elements $k \in [n]$ satisfying
$\bar{\delta}^m(-k) \neq k$. Because $\bar{\delta}^{-1}(-k) =
-\bar{\delta}(k)$ for all $k \in [n]$,
 there exist at most two $k$'s (one $k$'s, resp.) satisfying $\bar{\delta}^m(-k) = k$  in any orbit $O$ under $\langle \bar{\delta}\rangle$
 if the size $|O|$ is
even(odd, resp.). Since $|O(1)| = \bar{n}$, there exist at least
$\bar{n} - 2$ elements $k \in O(1)$ satisfying $\bar{\delta}^m(-k)
\neq k$. If there exists an orbit $O$ under $\langle
\bar{\delta}\rangle$ which is not $O(1)$ and whose size is greater
than or equal to 3 then there exists at least two $k$ in $O$ such
that $\bar{\delta}^m(-k) \neq k$. It implies that there exist at
least $\bar{n}$ elements $k \in [n]$ satisfying
$\bar{\delta}^m(-k) \neq k$, a contradiction. Therefore, except
$O(1)$, the size of every orbit under $\langle
\bar{\delta}\rangle$ is 1 or 2. By our assumption that there is an
orbit under $\langle \bar{\delta}\rangle$ whose size doesn't
divide $\bar{n}$, the value $\bar{n}$ should be
 odd. It means that there exist at least  $\bar{n}-1$ elements $k$ in $O(1)$ satisfying
$\bar{\delta}^m(-k) \neq k$. It implies that  there exist exactly $\bar{n} +
1$ elements $k \in [n]$ satisfying $\bar{\delta}^m(-k) = k$ and
$f(\bar{\delta})=1$, namely, $\Fix(\bar{\delta}) = \{0 \}$. Hence,
each orbit under $\langle \bar{\delta}\rangle$ containing neither
$0$ nor $1$ is $\{i,-i\}$ for some $i \in [n]$ and $m$ is odd.
Since $\bar{\delta}^{-1}(-k) = -\bar{\delta}(k)$ for all $k \in
[n]$ and there exists one $k \in O(1)$ such that
$\bar{\delta}^m(-k) = k$, it holds that
  $-O(1) = \{ -k \ | \
k \in O(1) \} = O(1)$. Recall that for any $k \in [n]$,
$\bar{\delta}^{\bar{n}}(k)=k$ if and only if $k \in O(1) \cup
\{ 0 \}$. For any orbit $\{i,-i\}$ under $\langle
\bar{\delta}\rangle$, $\bar{\delta}^{\bar{n}}(i+1)  =
\bar{\delta}^m(-i)+1 = i+1$ and $\bar{\delta}^{\bar{n}}(-i+1) =
\bar{\delta}^m(i)+1 = -i+1$.  It implies that $i-1,
i+1, -i-1, -i+1 \in O(1)$ because $i\neq \pm 1$ and $O(1) =
-O(1)$. Hence, there exist no two consecutive elements $i, i+1 \in
[n]$ satisfying $|O(i)|=|O(i+1)|=2$. Note that $|O(0)|=1$ and
$|O(1)|=|O(-1)|=\bar{n}$. Since there are $\bar{n}-1$ elements $j
\in [n]$ such that $|O(j)|=2$, for any even $2k \in [n]$,
$O(2k)=\{ 2k, -2k \}$ or equivalently $\bar{\delta}(2k)=-2k$.
Moreover, $O(1)$ is composed of all odd numbers.  Hence, for any
even $2k \in [n]$,
$$\bar{\delta}^m(2k+1) =\bar{\delta}^m(-(-2k-1)) =
\bar{\delta}^{\bar{n}}(-2k) -1 = 2k-1.$$ It implies that $m$ and
$\bar{n}$ are relative prime. Moreover, $m$ and $n$ are relative
prime because $m$ is odd. So there exists $s \in [n]$ such that
$sm \equiv 1 (\mod n)$. For any even $2k \in [n]$,
$$\bar{\delta}(2k+1) = \bar{\delta}^{sm}(2k+1)=\bar{\delta}^{(s-1)m}(2k-1) =
\cdots=2k+1-2s.$$ Let $x=-2s$. Then, $\bar{\delta} = (0)(2 \ -2)(4
\ -4)\cdots(1 \ \ 1+x \ \ 1+2x \ \ 1+3x \cdots)$.

By Lemma~\ref{disjoint-union}$(3)$, there exist $a(\bar{n})$ and
$b(\bar{n})$ satisfying Eq.$(\ast_1)$ or Eq.$(\ast_2)$ with
$\bar{\delta}$.
 Suppose that $a(\bar{n})$
and $b(\bar{n})$ satisfy Eq.$(\ast_1)$. Then, by
Lemma~\ref{a(i)b(i)},
 $$b(\bar{n}) = \bar{\delta}^{\bar{n}}(1) = 1 \ \ \mbox{and}
  \ \ b(\bar{n})=\bar{\delta}^{-\bar{\delta}(\bar{n})}(1) = \bar{\delta}^{-(\bar{n}+x)}(1) \neq 1,$$
  which is a contradiction. Hence, we can assume that $a(\bar{n})$
and $b(\bar{n})$ satisfy Eq.$(\ast_2)$.
  By Lemma~\ref{a(i)b(i)},
$$b(\bar{n}) =- \bar{\delta}^{\bar{n}}(1) = -1 \ \ \mbox{and}
  \ \ b(\bar{n})=\bar{\delta}^{-\bar{\delta}(\bar{n})}(1) = \bar{\delta}^{-(\bar{n}+x)}(1) = 1-(\bar{n}+x)x=1-x^2$$
 because $x$ is even.
 It implies that $x^2 \equiv 2 (\mod n)$. So, $\bar{\delta} \in \mathcal{N}^{non}_n$, which contradicts the assumption.
\end{pf}

 {\bf Remark}  As one can see in the proof of
Lemma~\ref{Lem=order-sigma}, any element in $\mathcal{M}^{non}_n$
which does not satisfy the condition in
Lemma~\ref{Lem=order-sigma} is $\bar{\delta}_{n,x}$ for some even
$n$ and $x$ such that $x^2 \equiv 2 (\mod n)$. This is a reason
why  we define $\bar{\delta}_{n,x}$  in Section 3.

\begin{prop}{\rm \cite{KK1}} \label{sigma-2n2} If $\bar{\delta}$ is the identity permutation of
$[n]$ then $|\left<R_{\bar{\delta}}, L \right> | = 2n^2$.
Furthermore, if we define $\bar{\delta} :[n] \rightarrow [n]$ by
 $\bar{\delta}(k) = k(1+rd)$ for all $k \in [n]$, where $n \geq 3$, $d$ is a divisor
of $n$ and  $r$ is a positive integer such that the order of
$1+rd$ in the multiplicative group $\mathbb{Z}_{n}^{*}$ of units is $d$, then
$|\left<R_{\bar{\delta}}, L \right> | = 2n^2$.
  \end{prop}
 \medskip

 \begin{lemma} \label{non-exi-d=2}
If $n \geq 3$, $
  |\langle \bar{\delta} \rangle | \neq 2$ for every $\bar{\delta} \in \mathcal{M}^{non}_n$.
 \end{lemma}

\begin{pf}
 Suppose that there
exists a $\bar{\delta} \in \mathcal{M}^{non}_n$ satisfying
  $|\langle \bar{\delta} \rangle |=2$. Then, $n$ is even.    By Lemma~\ref{disjoint-union}(3), there
exist $a(1), b(1) \in [n]$ satisfying Eq.$(\ast_1)$ or
Eq.$(\ast_2)$ with $\bar{\delta}$.  In both cases, $a(1)=\bar{\delta}(1)$ and
$b(1)=\bar{\delta}^{-\bar{\delta}(1)}(1)$ by Lemma~\ref{a(i)b(i)}.
 Suppose $a(1)=\bar{\delta}(1)$ is even and let $\bar{\delta}(1) = 2r$. Then, $b(1)=1$ and $$ \bar{\delta}(k + 1) =
 \bar{\delta}^{b(1)}(k) + a(1)=\bar{\delta}(k) + 2r \ \ \mbox{for all $k \in [n]$ \ \ \ or} $$  $$\bar{\delta}(k + 1) =
 \bar{\delta}^{b(1)}(-k) + a(1)=\bar{\delta}(-k) + 2r = -\bar{\delta}(k) + 2r \ \ \mbox{for all $k \in [n]$}.$$
In both cases, one can inductively show that $\bar{\delta}(k)$ is
even for all $k \in [n]$. It contradicts that $\bar{\delta} \in
S_{[n]}$. Therefore, we can assume that $a(1)=\bar{\delta}(1)$ is
odd. Let $\bar{\delta}(1) = 1+ 2r$. Then,
$b(1)=\bar{\delta}^{-\bar{\delta}(1)}(1)=\bar{\delta}(1)=1+2r$ by
Lemma~\ref{a(i)b(i)}.

 Suppose that Eq.$(\ast_1)$ holds.
 Then,
  $$\bar{\delta}(k + 1) = \bar{\delta}^{b(1)}(k) + a(1) = \bar{\delta}(k) +
 1+2r = \bar{\delta}(k-1) + 2(1+2r)
 =\cdots = (k+1)(1+2r).$$
 Moreover, $2$ is the smallest positive integer $d$ satisfying
 $\bar{\delta}^d(1)=(1+2r)^d = 1$.
 By Proposition~\ref{sigma-2n2}, $|\left<R_{\bar{\delta}}, L \right> | =
 2n^2$. So, $\bar{\delta} \notin \mathcal{M}^{non}_n$, a contradiction.

 Now, suppose that Eq.$(\ast_2)$ holds.
 Then, $b(1)=- \bar{\delta}(1)$ and hence $b(1)=\bar{\delta}(1)=- \bar{\delta}(1)$.  Since
 $-\bar{\delta}(1) = \bar{\delta}^{-1}(-1) = \bar{\delta}(-1)$, we have
 $\bar{\delta}(1)=\bar{\delta}(-1)$. It implies that  $n$ is $2$, a
 contradiction.
\end{pf}

 From
 Proposition~\ref{sigma-2n2} and Lemma~\ref{non-exi-d=2}, one can see that for any  $n \geq
3$ and for every $\bar{\delta} \in \mathcal{M}^{non}_n,$
$\bar{\delta}$ is neither the identity nor an involution.

\begin{lemma} \label{congruence}
 Suppose that $\bar{\delta} \in \mathcal{M}^{non}_n - \mathcal{N}^{non}_n$ with $n \geq
 3$ exists and let $|\langle \bar{\delta} \rangle|=d$. If $k_1 \equiv k_2 (\mod
 d)$ for some $k_1,k_2 \in [n]$ then $\bar{\delta}(k_1) \equiv \bar{\delta}(k_2) (\mod
 d)$.
\end{lemma}
\begin{pf}
 By Lemma~\ref{Lem=order-sigma}, $d$ is a divisor of $n$ and $d < n$.
By Lemma~\ref{disjoint-union}$(3)$, there exist $a(d)$ and $b(d)$
satisfying Eq.$(\ast_1)$ or Eq.$(\ast_2)$. Assume that
Eq.$(\ast_1)$ holds. Then, $b(d)=\bar{\delta}^{d}(1)=1$  by
Lemma~\ref{a(i)b(i)},
 which means $k+1 = \bar{\delta}^{d}(k)+1 = \bar{\delta}^{a(d)}(k + b(d)) = \bar{\delta}^{a(d)}(k + 1)$.
 It implies that $a(d)$ is a multiple of $d$, say $a(d) =
 rd$. So, the first equation in Eq.$(\ast_1)$ is  $\bar{\delta}(k + d) =
 \bar{\delta}^{b(d)}(k) + a(d) =
 \bar{\delta}(k) + rd$. Therefore, if $k_1 \equiv k_2 (\mod
 d)$ for some $k_1,k_2 \in [n]$, then $\bar{\delta}(k_1) \equiv \bar{\delta}(k_2) (\mod
 d)$.

 Next, suppose that Eq.$(\ast_2)$
holds. Then, $b(d)=-\bar{\delta}^{d}(1)=-1$  by
Lemma~\ref{a(i)b(i)}, which means $k+1 = \bar{\delta}^{d}(k)+1
 = \bar{\delta}^{a(d)}(-k + b(d)) = \bar{\delta}^{a(d)}(-k - 1)$.
By taking $k=-2$ and $k=-\bar{\delta}(1)-1$ in the equation
$\bar{\delta}^{a(d)}(-k - 1) =
 k+1$, one can get $\bar{\delta}^{a(d)}(1)=-1$ and $\bar{\delta}^{a(d)+1}(1) =-\bar{\delta}(1)$. Since
 $-\bar{\delta}(1) = \bar{\delta}^{a(d)+1}(1)=
 \bar{\delta}(\bar{\delta}^{a(d)}(1)) = \bar{\delta}(-1) = -\bar{\delta}^{-1}(1)$, we have $\bar{\delta}^{-1}(1) = \bar{\delta}(1)$. By Lemma~\ref{Lem=order-sigma},
 $\bar{\delta}^{-1} = \bar{\delta}$, or equivalently $d=1$ or $2$.  It is impossible by Proposition~\ref{sigma-2n2} and
 Lemma~\ref{non-exi-d=2}.
 \end{pf}

Suppose that $\bar{\delta} \in \mathcal{M}^{non}_n -
\mathcal{N}^{non}_n$ with $|\left< \bar{\delta}\right>| = d$
exists. By Lemma~\ref{congruence}, the function
$\bar{\delta}_{(1)} : [d] \rightarrow [d]$ defined by
$\bar{\delta}_{(1)} (k) \equiv \bar{\delta}(k) (\mod d)$ for any
$k \in [d]$ is well-defined. Furthermore, $\bar{\delta}_{(1)}$ is
a bijection, namely, a permutation of $[d]$. We call the
permutation $\bar{\delta}_{(1)}$ the $(\mod d)$-\emph{reduction}
of $\bar{\delta}$.
  In fact, $\bar{\delta}_{(1)}$ belongs to
$\mathcal{M}^{non}_d$ as the following lemma in a general setting.

\begin{lemma} \label{sigma'alsoMn-non}
 Suppose that $\bar{\delta} \in \mathcal{M}^{non}_n -
\mathcal{N}^{non}_n$ with $|\left< \bar{\delta}\right>| = d \geq
3$ exists. Let $m$ be a divisor of $n$ such that \vspace{-.2cm}
\begin{enumerate} \item[$(1)$]  $m$ is a multiple of $d$ and \vspace{-.2cm}
\item[$(2)$] if $k_1 \equiv k_2 (\mod
 m)$  for some $k_1,k_2 \in [n]$ then $\bar{\delta}(k_1) \equiv
\bar{\delta}(k_2) (\mod
 m)$.
 \end{enumerate} Define $\bar{\delta}'
:[m] \rightarrow [m]$ by $\bar{\delta}' (k) \equiv \bar{\delta}(k)
(\mod m)$ for any $k \in [m]$. Then, $\bar{\delta}'$ is a
well-defined bijection and it belongs to $\mathcal{M}^{non}_m$.
\end{lemma}

\begin{pf} By the assumption that $\bar{\delta}(k_1) \equiv \bar{\delta}(k_2) (\mod
 m)$ for any $k_1,k_2 \in
[n]$ satisfying $k_1 \equiv k_2 (\mod
 m)$, $\bar{\delta}'
:[m] \rightarrow [m]$ is well-defined. Since  $\bar{\delta}$ is a
bijection,  $\bar{\delta}'$ is also a bijection. By the fact
$\bar{\delta}(0)=0$, we have $\bar{\delta}'(0) = 0$. It is easily
checked that $(\bar{\delta}')^{-1}(m-k)=m-\bar{\delta}'(k)$ for
any $k \in [m]$.

Now, we aim to show that $\bar{\delta}' \in \mathcal{M}^{non}_m$
using Lemma~\ref{disjoint-union}$(3)$.  For any $k \in [n]$, let
$k'$ denote the remainder of $k$ divided by $m$. By Lemma
\ref{disjoint-union}(3), for any $i \in [n]$ there  exist $a(i)$
and $b(i)$ satisfying Eq.$(\ast_1)$ or Eq.$(\ast_2)$. One can
easily show that if we define $a(i')=a(i)'$ and $b(i')=b(i)'$ then
$a(i')$ and $b(i')$ also satisfy Eq.$(\ast_1)$ or Eq.$(\ast_2)$
depending on whether $a(i)$ and $b(i)$ satisfy Eq.$(\ast_1)$ or
Eq.$(\ast_2)$.   Since $\bar{\delta} \in
 \mathcal{M}^{non}_n$, there exists at least one $j \in [n]$ such that  $\bar{\delta}(k + j) =
 \bar{\delta}^{b(j)}(-k) + a(j)$ and $\bar{\delta}^{j}(k) +1 =
 \bar{\delta}^{a(j)}(-k + b(j))$ for all $k \in [n]$ by
 Lemma~\ref{disjoint-union}$(3)$. It implies that $\bar{\delta}'(k' + j') \equiv
 \bar{\delta}'^{b(j)'}((-k)') + a(j)' (\mod m)$ and $\bar{\delta}'^{j'}(k') +1 \equiv
 \bar{\delta}'^{a(j)'}((-k)' + b(j)') (\mod m)$ for all $k' \in [m]$. So, by
 Lemma~\ref{disjoint-union}, $\bar{\delta}' \in
 \mathcal{M}^{non}_m$.
\end{pf}

\begin{corollary} \label{sigma(1)alsoMn-non}
Suppose that $\bar{\delta} \in \mathcal{M}^{non}_n -
\mathcal{N}^{non}_n$ with $|\left< \bar{\delta}\right>| = d \geq
3$ exists. Then, $\bar{\delta}_{(1)}$ belongs to
$\mathcal{M}^{non}_d$.
\end{corollary}
\begin{pf}
By Lemmas~\ref{congruence} and \ref{sigma'alsoMn-non},  the $(\mod
d)$-reduction $\bar{\delta}_{(1)}$ of $\bar{\delta}$ belongs to
$\mathcal{M}^{non}_d$
\end{pf}

\section{Proof of Theorem~\ref{main-theorem}}

To prove Theorem~\ref{main-theorem}, we need to show that for any
integer $n \equiv 0, 1$ or $3 (\mod 4)$,  no nonorientable regular
embedding of $K_{n,n}$ exists and for $n \equiv 2(\mod 4)$,
$\mathcal{M}^{non}_n = \mathcal{N}^{non}_n$.

For a non-negative integer $k$,   we define $\bar{\delta}_{(0)} =
\bar{\delta}\in S_{[n]}$ and $\bar{\delta}_{(k+1)} =
(\bar{\delta}_{(k)})_{(1)}$ by taking reduction inductively.

\begin{lemma} \label{non-exi-sig(1)id}
Suppose that $\bar{\delta} \in \mathcal{M}^{non}_n
-\mathcal{N}^{non}_n$ with $n \geq 3$  would exist. Then,
\begin{enumerate}
 \item[$(1)$] $\bar{\delta}_{(1)}$ is not the identity, and
  \item[$(2)$] $|\langle \bar{\delta}
\rangle|$ is even.
 \end{enumerate}
 \end{lemma}

\begin{pf}  Suppose that there
exists a $\bar{\delta} \in \mathcal{M}^{non}_n
-\mathcal{N}^{non}_n$.  Let $|\langle \bar{\delta} \rangle |=d$.  By
Proposition~\ref{sigma-2n2} and Lemma~\ref{non-exi-d=2}, $d \geq 3$. It implies that $\bar{\delta}_{(1)}$ belongs to
$\mathcal{M}^{non}_d$ by Corollary~\ref{sigma(1)alsoMn-non}. Hence,  $\bar{\delta}_{(1)}$ is not the identity  by Proposition~\ref{sigma-2n2}.

 Suppose that $|\langle
\bar{\delta} \rangle |=d$ is odd.    By
Lemma~\ref{Lem=order-sigma}, $d$ is less than $n$. Since for any
odd $n$, $\mathcal{N}^{non}_n = \emptyset$, $\bar{\delta}_{(1)}$
is an element in $\mathcal{M}^{non}_d -\mathcal{N}^{non}_d$ by
Corollary~\ref{sigma(1)alsoMn-non} and the order of
$\bar{\delta}_{(1)}$ is also odd. By continuing the same process,
one can get  $j \geq 1$ and $d_j \geq 3$ such that
$\bar{\delta}_{(j)} \in \mathcal{M}^{non}_{d_{j}}
-\mathcal{N}^{non}_{d_{j}}$ and
 $\bar{\delta}_{(j+1)}$ is the identity permutation on $[d_{j+1}]$,
where $d_{j+1} = |\left< \bar{\delta}_{(j)} \right>|$ and $d_{j} =
|\left< \bar{\delta}_{(j-1)} \right>| \geq 3$. But, such
$\bar{\delta}_{(j)}$ cannot exist by (1).
\end{pf}

\begin{corollary} \label{no-emb-for-odd}
 If $n$ is odd, $\mathcal{M}^{non}_n = \emptyset$, or
equivalently there is no nonorientable regular embedding of $K_{n,n}$.
\end{corollary}
\begin{pf}
Suppose that $\bar{\delta} \in \mathcal{M}^{non}_n$ exists. Since
$\mathcal{N}^{non}_n = \emptyset$ for odd $n$, $\bar{\delta}$
belongs to $\mathcal{M}^{non}_n - \mathcal{N}^{non}_n$. By
Lemma~\ref{Lem=order-sigma}, the order $|\langle \bar{\delta}
\rangle|$ is a divisor of $n$. Hence, $|\langle \bar{\delta}
\rangle|$ is odd, which is a contradiction by
Lemma~\ref{non-exi-sig(1)id}.
\end{pf}

\begin{lemma} \label{non-exi-sig(1)in-nmp}
There does not exist $\bar{\delta} \in \mathcal{M}^{non}_n -
\mathcal{N}^{non}_n$ with $|\left< \bar{\delta} \right> | = d \geq
3$ such that $\bar{\delta}_{(1)} \in \mathcal{N}^{non}_d$.
 \end{lemma}
 \begin{pf}
 Suppose that there exists an element $\bar{\delta} \in \mathcal{M}^{non}_n - \mathcal{N}^{non}_n$ of order $d \geq 3$ such that $\bar{\delta}_{(1)} \in
\mathcal{N}^{non}_d$.  Note that $d \equiv 2 (\mod 4)$. We
consider two cases  that $n \equiv 2 (\mod 4)$ and $n \equiv 0
(\mod 4)$ separately.
\medskip

\noindent \emph{Case 1. $n \equiv 2 (\mod 4)$.} \\
Then, $\bar{n}$ is an odd integer. Let $O$ be the orbit of
$\bar{n}$ under $\langle \bar{\delta}\rangle$. Then, the size
$|O|$ is a divisor of $d=|\left< \bar{\delta} \right> |$.
Furthermore, $|O|$ is a multiple of $d/2$ because all odd numbers
in $[d]$ are in the same orbit under $\langle \bar{\delta}_{(1)}
\rangle$ whose size is  $d/2$. Therefore, $|O|$ is $d/2$ or $d$.
Since $-\bar{n} = \bar{n}$ and $\bar{\delta}^{-1}(-k) =
-\bar{\delta}(k)$ for any $k \in [n]$, one can see that $-O = O$
and the size $|O|$ is odd, which implies $|O| =  d/2$. Since all
odd numbers in $[d]$ are in the same orbit under $\langle
\bar{\delta}_{(1)} \rangle$, there exists a number $1+jd \in [n]$
such that $1+jd \in O$. It implies that $\bar{\delta}^{d/2}(1+jd)
= 1+jd$ and hence $(L R_{\bar{\delta}}^{1+jd} L)^{-1}
R_{\bar{\delta}}^{d/2} (L R_{\bar{\delta}}^{1+jd} L) (0) = 0$. So,
as a conjugate of $R_{\bar{\delta}}^{ d/2}$, $(L
R_{\bar{\delta}}^{1+jd} L)^{-1} R_{\bar{\delta}}^{d/2} (L
R_{\bar{\delta}}^{1+jd} L)$ belongs to the vertex stabilizer $
\left<R_{\bar{\delta}}, L \right>_{0} = \left<R_{\bar{\delta}}, t
\right>$ which is isomorphic to dihedral group $D_n$ of order
$2n$. Since the order of $(L R_{\bar{\delta}}^{1+jd} L)^{-1}
R_{\bar{\delta}}^{ d/2} (L R_{\bar{\delta}}^{1+jd} L)$ is not 2,
$(L R_{\bar{\delta}}^{1+jd} L)^{-1} R_{\bar{\delta}}^{d/2} (L
R_{\bar{\delta}}^{1+jd} L) = R_{\bar{\delta}}^{m}$ for some $m \in
[n]$. Because $R_{\bar{\delta}}^{m}$ and $R_{\bar{\delta}}^{ d/2}$
are conjugate in $\left<R_{\bar{\delta}}, L\right>$, they have the
same order and consequently, $\langle R_{\bar{\delta}}^{m} \rangle
=\langle R_{\bar{\delta}}^{ d/2} \rangle$ as subgroups of the
cyclic group $\left<R_{\bar{\delta}}\right>$. Since $d/2$ is a
divisor of $n$, there exists $ \ell \in [n/\frac{d}{2}]$ such that
$m = \ell  d/2$ and $( \ell , n/\frac{d}{2}) = 1$. By considering
two images of $1+jd$ under the permutations $(L
R_{\bar{\delta}}^{1+jd} L)^{-1} R_{\bar{\delta}}^{ d/2} (L
R_{\bar{\delta}}^{1+jd} L)$ and $R_{\bar{\delta}}^{\ell
 d/2}$, we have
$$\bar{\delta}^{ d/2}(2+2jd) -1-jd = (L R_{\bar{\delta}}^{1+jd} L)^{-1}
R_{\bar{\delta}}^{ d/2} (L R_{\bar{\delta}}^{1+jd} L)(1+jd) =
R_{\bar{\delta}}^{\ell  d/2}(1+jd) = 1+jd.$$ It implies that
$\bar{\delta}^{ d/2}(2+2jd)= 2+2jd$. Since $ d/2$ is odd and for
any even $k \in [d]$ with $k\neq 0$, the orbit of $k$ under
$\langle \bar{\delta}_{(1)} \rangle$
 is $\{k, d-k\}$, the even number $2+2jd$ should be a
multiple of $d$. It means that $d$ is 1 or 2, a contradiction.
Therefore, for any $n \equiv 2 (\mod 4)$, no $\bar{\delta} \in
\mathcal{M}^{non}_n - \mathcal{N}^{non}_n$ with $|\left<
\bar{\delta} \right> | = d \geq 3$ such that $\bar{\delta}_{(1)}
\in \mathcal{N}^{non}_d$ exists.

\bigskip

\noindent \emph{Case 2. $n \equiv 0 (\mod 4)$.} \\
Let $n = 2sd$ for some even integer $2s$. Then, $\bar{n}=sd$ is
even. By Lemma \ref{disjoint-union}(3), there  exist $a(\bar{n})$
and $b(\bar{n})$ satisfying Eq.$(\ast_1)$ or Eq.$(\ast_2)$.
 In both cases,
 $a(\bar{n})=\bar{\delta}(\bar{n})$ by
 Lemma~\ref{a(i)b(i)}. Since $\bar{n}$ is a multiple of $d$,
 $a(\bar{n})$ is also a multiple of $d$ by Lemma~\ref{congruence}.

 Suppose that Eq.$(\ast_2)$ holds.
Then,
 $$k+1 = \bar{\delta}^{\bar{n}}(k)+1
 = \bar{\delta}^{a(\bar{n})}(-k + b(\bar{n}))=-k+b(\bar{n})$$ for all $k \in
 [n]$. It means that $b(\bar{n})=2k+1$ for all $k \in
 [n]$. Since $b(\bar{n})$ is a constant, $n \leq 2$, a
 contradiction.  So, Eq.$(\ast_1)$ holds, that is, $\bar{\delta}(k + \bar{n}) =
 \bar{\delta}^{b(\bar{n})}(k) + a(\bar{n})$ and $k+1 = \bar{\delta}^{\bar{n}}(k)+1 = \bar{\delta}^{a(\bar{n})}(k + b(\bar{n}))
 =k+b(\bar{n})$ for all $k \in [n]$. It means that $b(\bar{n})=1$. Hence,
 $\bar{\delta}(k + \bar{n}) =
 \bar{\delta}^{b(\bar{n})}(k) + a(\bar{n}) =
 \bar{\delta}(k) + \bar{\delta}(\bar{n})$. By taking $k=\bar{n}$ in the
 above equation, we have $2\bar{\delta}(\bar{n})=0$. Since
 $\bar{\delta}(\bar{n}) \neq 0$,  $\bar{\delta}(\bar{n}) =
 \bar{n}$ and $$\bar{\delta}(k + \bar{n}) =
 \bar{\delta}(k) + \bar{\delta}(\bar{n})= \bar{\delta}(k) + \bar{n}.$$ It
 implies that if $k_1 \equiv k_2 (\mod
 \bar{n})$ then $\bar{\delta}(k_1) \equiv \bar{\delta}(k_2) (\mod
 \bar{n})$. Let $\bar{\delta}':[\bar{n}] \rightarrow
 [\bar{n}]$ be defined by $\bar{\delta}'(k) \equiv \bar{\delta}(k) (\mod \bar{n})$ for any $k \in
 [\bar{n}]$. Then, by Lemma~\ref{sigma'alsoMn-non}, $\bar{\delta}'$ is  well-defined and it belongs to
  $\mathcal{M}^{non}_{\bar{n}}$ because $\bar{n}$ is a multiple of $d$. Note that the size of the
 orbit of 1 under $\langle \bar{\delta}' \rangle$  is $ d/2$ or $d$.
\medskip

\noindent \emph{Subcase 2.1.  The size of the
 orbit of 1 under $\langle \bar{\delta}' \rangle$ is $ d/2$.}  \\
 Let $d' = |\langle \bar{\delta}' \rangle|$. Then, $d'$ is $ d/2$ or $d$.  Since $d$ is a divisor of
 $\bar{n}$ and the orbit of 2 under $\langle \bar{\delta}_{(1)} \rangle$ is $\{ 2, d-2
 \}$, the size of the orbit of 2 under $\langle \bar{\delta}' \rangle$
  is even. Hence, $d'$ is even and consequently equals to $d$.
  Since the
order of $\bar{\delta}'$ is not equal to the size of the orbit of
1 under $\langle \bar{\delta}' \rangle$, $\bar{\delta}' \in
\mathcal{N}^{non}_{\bar{n}}$ by Lemma~\ref{Lem=order-sigma}.
Hence, $\bar{\delta}' = \bar{\delta}_{\bar{n}, x}$ for some $x \in
[\bar{n}]$ satisfying $x^2 \equiv 2 (\mod \bar{n})$. Moreover, $
d'=d=\bar{n}$. It implies that
$\bar{\delta}' = \bar{\delta}_{(1)}$ and all odd numbers in $[n]$
belong to the same orbit under $\langle \bar{\delta} \rangle$.
Furthermore, for any even number $2k \in [n]\setminus \{0, \bar{n}
\}$,  the size of the orbit of $2k$ under $\langle \bar{\delta}
\rangle$ is 2 or 4. Since it is a divisor of $d$ and $d \equiv 2
(\mod 4)$, it is 2. Note that the orbit of $2k$ under $\langle
\bar{\delta} \rangle$ is $\{2k, -2k \}$ or $\{2k, \bar{n}-2k \}$.

First, we want to show that $\bar{\delta}(2k)=-2k$ for all even
number $2k \in [n]$. By Lemma~\ref{disjoint-union}$(3)$, there
exist $a(2)$ and $b(2)$ satisfying Eq.$(\ast_1)$ or Eq.$(\ast_2)$.
In both cases,
 $a(2) = \bar{\delta}(2)$ and $b(2)=\bar{\delta}^{-a(2)}(1)=\bar{\delta}^{-\bar{\delta}(2)}(1)=\bar{\delta}^{2}(1)$ by Lemma~\ref{a(i)b(i)}.

Suppose that Eq.$(\ast_2)$ holds.
By Lemma~\ref{a(i)b(i)}, $b(2)=-\bar{\delta}^{2}(1)$. Hence,
$b(2)=-\bar{\delta}^{2}(1)=\bar{\delta}^{2}(1)$, which means
$2\bar{\delta}^{2}(1)=0$. Since $\bar{\delta}^{2}(1)$ is not 0,
$\bar{\delta}^{2}(1) = \bar{n}$. It contradicts the fact that the
orbit of 1 under $\langle \bar{\delta} \rangle$ is composed of all
odd numbers in $[n]$. So, Eq.$(\ast_1)$ holds.
By Lemma~\ref{a(i)b(i)}, $b(2) = \bar{\delta}^2(1) \equiv 1+2x
(\mod \bar{n})$. So, the first equation in  Eq.$(\ast_1)$ can be
written by $\bar{\delta}(k + 2) =
 \bar{\delta}^{b(2)}(k) + a(2)= \bar{\delta}^{1+2x}(k) +\bar{\delta}(2)$. Suppose that $\bar{\delta}(2)=\bar{n}-2$. Then,  $$\bar{\delta}(k+2)=\bar{\delta}^{1+2x}(k)
 +\bar{n}-2.$$  Taking
 $k=2$ in the equation $\bar{\delta}(k+2)=\bar{\delta}^{1+2x}(k)
 +\bar{n}-2$, we have
$\bar{\delta}(4)=\bar{\delta}^{1+2x}(2) +\bar{n}-2=
\bar{\delta}(2)+\bar{n}-2 =
 -4$.
 Taking $k=4$, we have $\bar{\delta}(6)= \bar{\delta}(4)+\bar{n}-2 =
 \bar{n}-6$. By continuing the same process, one can see that
 $\bar{\delta}(4k)=-4k$ and $\bar{\delta}(4k+2) = \bar{n}-4k-2$. Since
 $\bar{n} \equiv 2 (\mod 4)$, we have
 $\bar{\delta}(\bar{n}) = \bar{n} - \bar{n} = 0$, which is a
 contradiction. Hence, $\bar{\delta}(2)=-2$, namely, it holds that $$\bar{\delta}(k+2)=\bar{\delta}^{1+2x}(k)
 -2.$$ Taking
 $k=2$ in the equation $\bar{\delta}(k+2)=\bar{\delta}^{1+2x}(k)
 -2$, we have $\bar{\delta}(4) =\bar{\delta}^{1+2x}(2)
 -2 = \bar{\delta}(2)-2 = -4$.
 Taking $k=4$, we have $\bar{\delta}(6)= \bar{\delta}(4)-2 =
 -6$. By continuing the same process, one can see that
 $\bar{\delta}(2k)=-2k$ for all even numbers $2k \in [n]$.

 Now, we aim to apply  Lemma~\ref{disjoint-union}$(3)$
  once more to show that Subcase 2.1 cannot happen.  There  exist $a(1)$ and
 $b(1)$ satisfying Eq.$(\ast_1)$ or Eq.$(\ast_2)$.
In both cases,  $a(1) = \bar{\delta}(1) \equiv 1+x (\mod
\bar{n})$.
 For our convenience, let $\bar{\delta}(1) = 1 + x_1$.

Suppose that Eq.$(\ast_1)$ holds. Then, $b(1)= \bar{\delta}(1) = 1
+ x_1$. So, it holds that $$\bar{\delta}(k + 1) =
 \bar{\delta}^{b(1)}(k) + a(1)= \bar{\delta}^{b(1)}(k) + \bar{\delta}(1)=\bar{\delta}^{1+x_1}(k)+1+x_1.$$  By
 taking $k=2$, we have
 $\bar{\delta}(3) = \bar{\delta}^{1+x_1}(2) +1+x_1 =  -2 + 1+x_1 \equiv x-1 (\mod \bar{n})$.
 Since $\bar{\delta}(3) \equiv 3+ x (\mod \bar{n})$, $4 \equiv 0 (\mod \bar{n})$.
  By the assumption that $\bar{n} \equiv 2 (\mod 4)$, $\bar{n} = d=2$, which contradicts the assumption
 that $d \geq 3$.
So, Eq.$(\ast_2)$ holds.
Hence, $b(1) = -\bar{\delta}(1) = -1-x_1$ and it holds that
$$\bar{\delta}(k)+1 = \bar{\delta}^{a(1)}(-k + b(1)) =
\bar{\delta}^{1+x_1}(-k-1-x_1)$$  for all $k \in [n]$.
 By
taking odd number $2k+1 \in [n]$, we have $\bar{\delta}(2k+1)+1 =
\bar{\delta}^{1+x_1}(-2k-2-x_1)=2k+2+x_1$. Hence,
$\bar{\delta}(2k+1)=2k+1+x_1$. By taking $k=2$ in  the equation
$\bar{\delta}(k)+1 = \bar{\delta}^{1+x_1}(-k-1-x_1)$, we get
$$-1 =
\bar{\delta}(2)+1 = \bar{\delta}^{1+x_1}(-3-x_1) = -3-x_1
+(1+x_1)x_1 = -3+x_1^2.$$ So, $x_1^2 = 2 (\mod n)$. It is
impossible because $n \equiv 0 (\mod 4)$.  Hence, Subcase 2.1
cannot happen.

\bigskip

\noindent \emph{Subcase 2.2. The size of the
 orbit of 1 under $\langle \bar{\delta}' \rangle$ is $d$.} \\
Since the order of $\bar{\delta}'$ divides that of $\bar{\delta}$,
the order of $\bar{\delta}'$ is $d$ which equals the size of the
 orbit of 1 under $\langle \bar{\delta}' \rangle$. It implies that
 $\bar{\delta}' \in \mathcal{M}^{non}_{\bar{n}} - \mathcal{N}^{non}_{\bar{n}}$.
  Moreover, since $d$ divides $\bar{n}$, it holds that $\bar{\delta}'_{(1)}=\bar{\delta}_{(1)} \in
\mathcal{N}^{non}_{d}$. Since the subcase 2.1 cannot happen, by
repeating
 the same process continually, one can get $n_1 \equiv 2 (\mod 4)$
 and $\tilde{\bar{\delta}} \in \mathcal{M}^{non}_{n_1} - \mathcal{N}^{non}_{n_1}$ with
 $|\langle
\tilde{\bar{\delta}} \rangle | = d =|\langle \bar{\delta}
\rangle|$ such that $\tilde{\bar{\delta}}_{(1)}
=\bar{\delta}_{(1)} \in \mathcal{N}^{non}_d$. But, it returns to
Case 1.
 \end{pf}

\bigskip

Now, we  prove Theorem~\ref{main-theorem}.  We know that there
exists only one nonorientable regular embedding of $K_{2,2}$ into
the projective plane. Let $n \geq 3$.

Suppose that $\mathcal{M}^{non}_n \varsupsetneq \mathcal{N}^{non}_n$
and let $\bar{\delta} \in \mathcal{M}^{non}_n - \mathcal{N}^{non}_n$ and  $|\left< \bar{\delta} \right>|=d$. Note that $d
< n$. By Lemma~\ref{non-exi-sig(1)id}  and Lemma~\ref{non-exi-d=2},
$\bar{\delta}_{(1)}$ is not the identity and $d \geq 3$ is even. By
Lemmas \ref{sigma(1)alsoMn-non} and \ref{non-exi-sig(1)in-nmp},
$\bar{\delta}_{(1)} \in \mathcal{M}^{non}_d -\mathcal{N}^{non}_d$.
 By continuing the same process, one can get  $j \geq 1$
and $d_j \geq 3$ such that $\bar{\delta}_{(j)} \in
\mathcal{M}^{non}_{d_{j}} -\mathcal{N}^{non}_{d_{j}}$ and
 $\bar{\delta}_{(j+1)}$ is the identity permutation on $[d_{j+1}]$,
where $d_{j+1} = |\left< \bar{\delta}_{(j)} \right>|$ and $d_{j} =
|\left< \bar{\delta}_{(j-1)} \right>| \geq 3$. But, this is
impossible by Lemma~\ref{non-exi-sig(1)id}.
 Therefore, for any $n \geq 3$,  $\mathcal{M}^{non}_n = \mathcal{N}^{non}_n$.
 It means that for any integer $n \equiv 0, 1$ or $3 (\mod 4)$, no
nonorientable regular embedding of $K_{n,n}$ exists and for $n
\equiv 2 (\mod 4)$, $\mathcal{M}^{non}_n = \mathcal{N}^{non}_n$.
Hence, by Lemma \ref{card-Mnmp}, for $n = 2p_1^{a_1}p_2^{a_2}\cdots
p_k^{a_k}$  (the prime decomposition of $n$), the number of
nonorientable regular embeddings of $K_{n,n}$ up to isomorphism is
$2^k$ if $p_i \equiv \pm 1  (\mod 8)$ for all $i=1,2,\ldots,k$; 0
otherwise. {\hfill \qed \newline
\medskip}

{\bf Remark}  For any $\bar{\delta}_{n,x} \in
\mathcal{N}^{non}_n$, the covalency(face size) of its derived
nonorientable regular map $\cal M$ is the order of
$LR_{\bar{\delta}_{n,x}}$ which is in fact
 8. Hence, the number of faces of the map $\cal M$ is
$n^2/4$. By the Euler formula, the supporting surface of
$\cal M$ is nonorientable surface with $(3n^2-8n+8)/4$ crosscaps.


\begin{thebibliography}{99}
 \bibitem{AG} W. W. Adams and L. J. Goldstein,
Introduction to Number Theory, Prentice Hall, 1976.
\bibitem{DJKNS1} S.F. Du, G.A. Jones, J.H. Kwak, R. Nedela and M. \v{S}koviera,
         Regular embeddings of $K_{n,n}$ where $n$ is a power of 2. I: Metacyclic case, Europ. J. Combinatorics,
         {\bf 28}(6), 2007, 1595-1609.
\bibitem{DJKNS2} S.F. Du, G.A. Jones, J.H. Kwak, R. Nedela and M. \v{S}koviera,
         Regular embeddings of $K_{n,n}$ where $n$ is a power of 2. II: Non-metacyclic case, submitted.
\bibitem{GNSS} A. Gardiner, R. Nedela, J. \u{S}ir\'{a}\u{n} and M. \u{S}koviera,
         Characterization of graphs which underlie regular maps on closed surfaces,
       J. London Math. Soc.  {\bf 59}, 1999, 100-108.
\bibitem{Hu} B. Huppert, Endliche Gruppen I, Springer-Verlag, 1967.
\bibitem{J}  G.A. Jones,
         Regular embeddings of complete bipartite graphs:
         classification and enumeration,
       submitted.
\bibitem{JNS1} G.A. Jones,  R. Nedela and M. \v{S}koviera,
         Complete bipartite graphs with a unique regular embedding, J. Combin. Theory Ser. B, {\bf 98}(2), 2008, 241-248.
\bibitem{JNS2} G.A. Jones,  R. Nedela and M. \v{S}koviera,
          Regular embeddings of $K_{n,n}$ where $n$ is an odd prime power, Europ. J. Combinatorics,
           {\bf 28}(6), 2007, 1863-1875.
\bibitem{KK}  J.H. Kwak and Y.S. Kwon,
         Classification of reflexible regular embeddings and self-Petrie dual regular embeddings  of  the complete bipartite graphs,
        Discrete Math., {\bf 308}(11), 2008, 2156-2166.
\bibitem{KK1}  J.H. Kwak and Y.S. Kwon,
         Regular orientable embeddings of complete bipartite graphs,
      J. Graph Theory {\bf 50}(2), 2005, 105-122.
 \bibitem{KN}  Y.S. Kwon and R. Nedela,
         Non-existence of  nonorientable regular embeddings of $n$-dimensional cubes,
       Discrete Math., {\bf 307}, 2007, 511-516.
 \bibitem{N} R. Nedela,  Regular maps - combinatorial objects relating
     different fields of mathematics, J. Korean Math. Soc.  {\bf 38}(5), 2001, 1069-1105.
 \bibitem{NSZ} R. Nedela, M. \u{S}koviera and A. Zlato\v{s},  Regular embeddings of complete bipartite graphs,
   Discrete Math. {\bf 258}, 2002, 379-381.
\bibitem{Wi} S.E. Wilson,  Cantankerous maps and rotary embeddings of $K_{n}$,
    J. Combin. Theory Ser. B  {\bf 47}, 1989, 262-273.
   \end{thebibliography}
 \end{document}